\documentclass[draft, 10pt, twoside]{article}
\usepackage{amsfonts}
\usepackage{stmaryrd}
\usepackage{amssymb}
\usepackage{euscript}
\usepackage{amsthm}
\usepackage{latexsym,amsmath,amscd}
\usepackage{mathrsfs}
\usepackage{graphicx}
 \usepackage{dsfont}
\usepackage{enumerate}
\usepackage[normalem]{ulem}

\usepackage{color}
\usepackage{lineno}

\RequirePackage{fix-cm}
\makeatletter
\newcommand*\bigcdot{\mathpalette\bigcdot@{.5}}
\newcommand*\bigcdot@[2]{\mathbin{\vcenter{\hbox{\scalebox{#2}{$\m@th#1\bullet$}}}}}
\makeatother

\usepackage{fancyhdr}
\title{{\Large \bf Characteristion of honest times and semimartingales of class-$(\Sigma)$}}

\usepackage[colorlinks=false,hidelinks]{hyperref}
\usepackage[capitalise]{cleveref}

\newtheorem{thm}{Theorem}[section]

\newtheorem{lem}[thm]{Lemma}
\newtheorem{pro}[thm]{Proposition}
\newtheorem{cor}[thm]{Corollary}
\newtheorem{rem}[thm]{Remark}
\newtheorem{rems}[thm]{Remarks}
\newtheorem{ex}[thm]{Example}
\newtheorem{defi}[thm]{Definition}
\newtheorem{hyp}[thm]{Assumptions}

\newcommand{\be}{\begin{equation}}
\newcommand{\ee}{\end{equation}}
\newcommand{\bde}{\begin{displaymath}}
\newcommand{\ede}{\end{displaymath}}
\newcommand{\beq}{\begin{eqnarray*}}
\newcommand{\eeq}{\end{eqnarray*}}
\newcommand{\beqa}{\begin{eqnarray}}
\newcommand{\eeqa}{\end{eqnarray}}
\newcommand{\bel }{\left\{\begin{array}{ll}}
\newcommand{\eel}{\cr \end{array} \right.}

\newcommand{\seq}[1]{{\lbrace #1 \rbrace}}

\def\wt{\widetilde}

\def\F{{\cal F}}

\def\ff{{\mathbb F}}

\def\rr{{\mathbb R}}

\def\nn{{\mathbb N}}

\def\P{{\mathbb P}}

\def\I{\mathds{1}}

\newcommand{\bt}{\begin{thm}}
\newcommand{\et}{\end{thm}}
\newcommand{\bl}{\begin{lem}}
\newcommand{\el}{\end{lem}}
\newcommand{\bp}{\begin{pro}}
\newcommand{\ep}{\end{pro}}
\newcommand{\bcor}{\begin{cor}}
\newcommand{\ecor}{\end{cor}}

\newcommand{\bd}{\begin{defi} \rm }
\newcommand{\ed}{\end{defi}}
\newcommand{\brem }{\begin{rem} \rm }
\newcommand{\erem }{\end{rem}}
\newcommand{\brems }{\begin{rems} \rm }
\newcommand{\erems }{\end{rems}}
\newcommand{\bhyp }{\begin{hyp} \rm }
\newcommand{\ehyp }{\end{hyp}}
\newcommand{\bex}{\begin{ex} \rm }
\newcommand{\eex}{\end{ex}}

\newcommand{\llb}{\llbracket\,}
\newcommand{\rrb}{\,\rrbracket}

\author{Libo Li\\
\\ School of Mathematics and Statistics
\\ University of New South Wales
\\ NSW 2052, Australia}

\date{\today}

\begin{document}

\title{{\Large \bf Characterisation of honest times and optional semimartingales of class-$(\Sigma)$}}

\author{Libo Li\\
\\ School of Mathematics and Statistics
\\ University of New South Wales, Sydney
\\ NSW 2052, Australia}

\maketitle
\begin{abstract}
Given a finite honest time, we first show that the associated Az\'ema optional supermartingale can be expressed as the drawdown and the relative drawdown of some local optional supermartingales with continuous running supremum. The relative drawdown representation then allows us to provide a characterisation of finite honest times using a family of non-negative local optional supermartingales with continuous running supremum which converges to zero at infinity. Then we extend the notion of semimartingales of class-$(\Sigma)$ by allowing for jumps in its finite variation part of the semimartingale decomposition. This enables one to establish the Madan-Roynette-Yor option pricing formula for a larger class of processes, and finally, we apply the extended formula to the construction of finite honest times.

\end{abstract}
\vfill\break 
\tableofcontents
\vfill \eject

\section*{Introduction}
In this paper, we present in a uniform fashion some new results from two closely related topics. The first is on the representation of the additive and multiplicative decomposition of the Az\'ema supermartingale associated with finite honest times or last passage times, see \cref{honest}. The second is on semimartingales of class-$(\Sigma)$ and the Madan-Roynette-Yor option pricing formula.  For more applications of honest times and semimartingale of class-$(\Sigma)$ in mathematical finance, we refer interested readers to Nikeghbali and Platen \cite{NP}.

To motivate the study, we recall Nikeghbali and Yor \cite{NY2} have shown, under the assumptions that all martingales are continuous and that the given finite honest time $\tau$ avoids all stopping times, see \cref{A}, the following additive and multiplicative representations of the {\it Az\'ema supermartingale} associate with $\tau$ given by $Z_t:=\P(\tau > t\,|\,\F_t)$ holds, that is
\begin{align}
Z_t &= 1+ m_t - \sup_{s\leq t}m_s\label{eqa}\\
Z_t & = \frac{M_t}{\sup_{s\leq t} M_s}\label{eqm}
\end{align}
where $m$ is a continuous local martingale and $M$ is a non-negative continuous local martingale with the property that $\lim_{t\rightarrow \infty} M_t = 0$. In other words, the process $1- Z$ can be expressed as the {\it drawdown}\footnote{The drawdown of a process is the running supremum of the process minus the process itself.} of a local martingale $m$ and the {\it relative drawdown} of a non-negative local martingale $M$. Conversely, given a non-negative continuous local martingale $M$ such that $\lim_{t\rightarrow \infty} M_t = 0$, the Az\'ema supermartingale of the finite honest time $\tau : = \sup\{s: M_s = \sup_{u\leq s} M_u\}$ is of the form given in \eqref{eqm}. In applications, the multiplicative decomposition and representation of the Az\'ema supermartingale has recently received interest from credit risk modeling and the study of asymmetric information, e.g. Aksamit et al. \cite{ACDJ}, Fontana et al. \cite{FJS}, Zwierz \cite{Z} and Kardaras \cite{K2}.

The above representations of $Z$ and the corresponding characterisation of honest times through non-negative local martingales was again recovered in Nikeghbali and Platen \cite{NP} under only the assumption that the finite honest time $\tau$ avoids all stopping times. Under similar assumptions, the multiplicative representation was also studied in Kardaras \cite{K} and Acciaio and Penner \cite{BI}. To illustrate the extent of their results, a counterexample from \cite{ACDJ} was given in \cite{BI} to show that there exist finite honest times for which $1-Z$ cannot be expressed as the relative drawdown of a non-negative local martingale with continuous supremum, i.e. the representation \eqref{eqm} does not hold. This observation then led to Song \cite{S} where, for an arbitrary random time, necessary and sufficient conditions for the representation of the form \eqref{eqm} to hold was obtained. Here, to better illustrate that there exist finite honest time for which \eqref{eqm} does not hold, we also provide two simple counterexamples in \cref{e01} and \cref{exp2.2}.

In view of the counterexamples, the most important contribution of the paper is that we remove the last standing assumption that $\tau$ avoids all stopping times to fill the final gap in the literature on the existence and uniqueness of the additive and the multiplicative representations of the Az\'ema supermartingale associated with a finite honest time, which is then used to provide a complete characterisation of finite honest times.  More precisely, by combining \cref{t1}, \cref{p2.18} and \cref{p2.21}, we show that given an arbitrary finite honest time, instead of local martingales in \eqref{eqa} and \eqref{eqm}, the process $1-Z$ can be uniquely expressed as the drawdown of some local supermartingale and relative drawdown of some non-negative local supermartingale with continuous running supremum. In hindsight, the main obstacle in removing the assumption that $\tau$ avoids all stopping times is that one was too focused on $Z$ and have insisted that $m$ and $M$ should be local martingales. In fact, instead of the supermartingale $Z$, it is more natural to consider the {\it Az\'ema optional supermartingale} $\wt Z_t :=\P(\tau \geq t\, |\,\F_t) $ as any finite honest time $\tau$ can be expressed as the end of the optional set $\seq{\wt Z = 1}$ and the representations for $Z$ can be obtained by noticing that $Z = \wt Z_+$. The switch from $Z$ to $\wt Z$ is crucial, and by doing so, we can remove the assumption that $\tau$ avoids all stopping times and obtain representations of $\wt Z$ in the form given in \eqref{eqa} and \eqref{eqm}, with the key difference being that the local martingales $m$ and $M$ are replaced by {\it local optional supermartingales} which exhibits l\`agl\`ad trajectories.
The main technical difficulty faced in this study is that the process $\wt Z$ is in general not a c\`adl\`ag process. Therefore, the standard c\`adl\`ag semimartingale calculus cannot be applied and the techniques employed in \cite{BI,K,NP, NY2}, such as the Doob maximal identity and the Skorokhod reflection lemma, e.g. Lemma 2.1 and Lemma 2.4 in \cite{NY2}, for c\`adl\`ag functions are not directly applicable and one needs to seek alternative methods. 

The second topic considered here is semimartingales of class-($\Sigma$). The notion of class-($\Sigma$) was first introduced for positive continuous submartingales in Yor \cite{Y1} and later extended in \cite{CNP, N1,NY2, Y1,Y2} to semimartingales and more recently examined in Eyi-Obiang et al. \cite{EOM, EOMT} in the context of signed measures. Notably, the authors of \cite{CNP} have shown that the introduction of class-$(\Sigma)$ allows for a martingale proof of the Madan-Roynette-Yor formula, see e.g.\cite{MRY}, which established a link between the last passage time of zero of a semimartingale of class-$(\Sigma)$ and the price of a European option. 

In the current definition of semimartingales of class-$(\Sigma)$, the predictable process of finite variation in the semimartingale decomposition is continuous. In the context of the Az\'ema supermartingale associated with a finite honest time, this continuity assumption is equivalent to the assumption that the honest time avoids all stopping times which we have previously removed. Hence the main contributions in the second part of the paper are (i) we extend the notion of semimartingales of class-$(\Sigma)$ by allowing for jumps in the predictable process of finite variation part of the decomposition and (ii) we recover, under the extended definition, some existing results for semimartingales of class-$(\Sigma)$, in particular the Madan-Roynette-Yor formula, and apply them to the construction of honest times.

The structure of the paper is as follows. In \cref{tool} we introduce the necessary notations and tools for our study. In \cref{ad} we consider the existence and uniqueness of an additive and a multiplicative representation for $\wt Z$. To prove our main result in \cref{t1}, we first derive a multiplicative decomposition of $\wt Z$ in \cref{mdz} and then identify the required representation  in  \cref{t1}.
Our approach is inspired by the works of Az\'ema, Meyer and Y\oe urp in \cite{A2,M, M0, YM} on the multiplicative decomposition of positive submartingales. Also, we rely heavily on the finer properties of honest times exposed in Jeulin \cite{J2} and stochastic calculus for l\`agl\`ad semimartingales under {\it the usual conditions}, which can be obtained from the $\underline{\underline{\mathtt{A}}}$-semimartingale calculus developed in Lenglart \cite{LE} or the stochastic calculus for optional semimartingales developed in Gal'\v cuk \cite{G1,G2,G3}. For notational convenience, we adopt here the framework of Gal'\v cuk.

Having obtained, in  \cref{t1}, the existence of an additive and a multiplicative representation of $\wt Z$ and identified the key properties of the local optional supermartingales involved. In \cref{sec2.1}, by using an extension of the Doob maximal identity obtained in \cref{dme}, we extend the existing characterisation of finite honest times which avoids all stopping times to all finite honest times in \cref{cor2.1} and study the uniqueness of the multiplicative representation in \cref{p2.18}. Finally, we provide a l\`agl\`ad extension of the Skorokhod reflection lemma in \cref{srl}, which we use in \cref{p2.21} to obtain the uniqueness of the additive representation.

In \cref{os}, under the extended definition of class-($\Sigma$) given in \cref{d3.2}, we generalise existing results on semimartingales of class-$(\Sigma)$ from \cite{CNP}. Firstly, we show in \cref{l3.1} and \cref{l3.3} that if $X$ and $Y$ are processes of class-($\Sigma$) then $X^+$, $X^-$, $|X|$ and $XY$ are again of class-$(\Sigma)$, and that any positive optional submartingale of class-$(\Sigma)$ can be represented as the drawdown of some optional supermartingale. Secondly, we recover the Madan-Roynette-Yor option pricing formula in \cref{t3.1} and \cref{c3.1}. Lastly, as an application of the results obtained in \cref{os}, we illustrate in \cref{max} and \cref{e3.1}, a method to construct examples of finite honest times for which the additive and multiplicative representation of the Az\'ema supermartingale can be retrieved from \cref{t1}, but not from the results of \cite{BI, K, NP, NY2}.

For the reader's convenience, we collect in the appendix some useful definitions and results from the theory of enlargement of filtrations and stochastic calculus for optional semimartingales.

\section{Notations and Terminologies}\label{tool}
We work on a filtered probability space $(\Omega,\mathcal A,\mathbb F,\mathbb P)$, where $\mathbb{F}:=(\F_t)_{t\geq 0}$ denotes a filtration satisfying the {\it usual conditions}, we set $\F_\infty := \bigvee_{t\geq 0} \F_t \subset \mathcal{A}$ and all martingales are c\`adl\`ag. The main tool used in this work is the stochastic calculus for optional semimartingales developed under the {\it unusual conditions} in Gal'\v cuk \cite{G1, G2, G3}. We stress that we do not make use of the full power of the calculus as $\ff$ is assumed to satisfy the {\it usual conditions} and all martingales are taken to be c\`adl\`ag. The case where $\ff$ does not satisfy the usual conditions can potentially be of interest, however, that would first require a complete study of finite honest times under the unusual conditions as our results in \cref{ad} relies on existing results for honest times which are all obtained under the usual conditions. On the other hand, the notion of optional semimartingale of class-$(\Sigma)$ and the corresponding results in \cref{os} can most likely be easily extended, under the unusual conditions, as the results of Gal'\v cuk already treats non-c\`adl\`ag martingales. However, we are reluctant to do so as there is a lack of concrete examples and applications.

Given a real-valued process $X$, as convention, we set $X_{0-} = 0$ and $X_\infty = \lim_{t\rightarrow \infty} X_t$ a.s, if it exists. The running supremum and infimum process of $X$ are denoted by $\overline{X}_t := \sup_{s\leq t} X_s$ and $\underline{X}_t := \inf_{s\leq t} X_s$. Given a c\`adl\`ag non-decreasing function $a$ on $\rr_+$, we say that the measure $da$ is carried on a set $G$ if $\int_{[0,\infty)}\I_{G^c}(s) da(s) = 0$. The support of $a$, that is the smallest closed set in $\rr_+$ which carries $da$, is given by $S(a) := \seq{t \geq 0: \forall \epsilon >0,\, a(t - \epsilon) < a(t + \epsilon) > 0}$ and the left support of $a$ is given by $S^g(a):= \seq{t \geq 0: \forall \epsilon >0, \,a(t - \epsilon) < a(t) > 0}$, see page 61, Chapter IV of Jeulin \cite{J2}.  We denote by $\mathcal{T}$ the set of all stopping times and for $0\leq s< t<\infty$, $\mathcal{T}_{[s,t]}$ the set of all stopping times $T$ such that $s\leq T\leq t$. A stochastic process $X$ is said to be of class-$(D)$ if the family $\seq{X_T\I_\seq{T<\infty}, T \in \mathcal{T}}$ is uniformly integrable, and it is said to be of class-$(DL)$ if for every $0<t <\infty$, the family $\seq{X_{T}, T \in \mathcal{T}_{[0,t]}}$ is uniformly integrable. For any integrable variation process $V$, we denote the $\ff$-optional (predictable) projection of $V$ by $^{o}V$ ($^{p}V$) and the $\ff$-dual optional (predictable) projection of $V$ by $V^{o}$ ($V^{p}$).
From Corollary 5.31 in He et al. \cite{HWY} the process $\,^oV - V^o$ is a uniformly integrable $\ff$-martingale starting at zero and $\,^o(\Delta V) = \Delta V^o$ holds.

Under the usual conditions, an optional martingale is a c\`adl\`ag uniformly integrable martingale, an optional local martingale is a c\`adl\`ag local martingale and any optional semimartingale $X$ takes the form $X = X_0 + M + A$, where $M$ is a c\`adl\`ag local martingale and $A$ is a l\`agl\`ad process of finite variation. As a convention, we suppose both $M$ and $A$ take value zero at time zero. We shall write $M^X$ and $A^X$ whenever there is a need to stress the dependence on $X$. In this setting, stochastic integrals for optional semimartingales reduces to the usual stochastic integrals, and one needs only to take care in counting the jumps of the integral against the process of finite variation $A$. As an alternative, one can apply the $\underline{\underline{\mathtt{A}}}$-semimartingale calculus in Lenglart \cite{LE} by taking $\underline{\underline{\mathtt{A}}} = \mathcal{O}(\ff)$, i.e. the optional $\sigma$-algebra generated by $\ff$, and the It\^o formula together with the solution to the stochastic exponential are readily available in section VI within. However, although it is more natural to apply the $\underline{\underline{\mathtt{A}}}$-semimartingale calculus as formulae developed in \cite{LE} are directly applicable under the usual conditions, we find the notations and presentations of Gal'\v cuk better suited for this work.

In the rest of this paper, unless otherwise stated, all stochastic processes in concern are optional semimartingales which are known to exhibit finite left and right limits. Given any l\`agl\`ad process $X$ we denote by $X_-$ and $X_+$ the left and right limits of $X$. The left and right jumps of $X$ are denoted by $\Delta X= X - X_-$ and $\Delta^+ X = X_+-X$ respectively. Any l\`agl\`ad process of finite variation $V$ can be decomposed into its right continuous part and left continuous part respectively by setting $V^g = \sum_{s< \cdot} \Delta^+ V_s$ and $V^r := V - V^g$. The right continuous part $V^r$ can be further decomposed into $V^r = V^c + V^d$, where $V^d = \sum_{s\leq \cdot} \Delta V_s$ and $V^c := V^r - V^d$. This gives us the decomposition
\begin{align}
V = V^c + V^d+ V^g. \label{Xg}
\end{align}
Finally, we mention that prior Gal'\v cuk \cite{G1} and Lenglart 	\cite{LE}, Mertens \cite{MJ} introduced under the usual conditions the notion of strong optional supermartingale and generalized the Doob-Meyer decomposition to this setting. However, here we work with optional supermartingales as defined by Gal'\v cuk and point out that the process $\wt Z$ is both an optional supermartingale and a strong optional supermartingale. Therefore, we will abuse slightly the terminology and call the Doob decomposition for optional supermartingales, the Doob-Meyer-Mertens-Gal'\v cuk decomposition. For more details on the general theory for stochastic processes, the reader is referred to He et al. \cite{HWY}, for results from the theory of enlargement of filtrations to Jeulin \cite{J2}. The reader can also refer to the recent book of Aksamit and Jeanblanc \cite{AJ} for a modern english exposition of the results from the theory of enlargement of filtration.

\section{The Additive and Multiplicative Representations}\label{ad}
In this part of the paper, given a finite honest time $\tau$, we study in \cref{s2.1}, the existence of an additive and multiplicative representation of the Az\'ema supermartingale associated with $\tau$ and, we study in \cref{sec2.1}, the uniqueness of such representations and provide a complete characterisation of finite honest time through a family of optional supermartingales.

\bd\label{honest}
A random time $\tau$ is a honest time, if for all $t\geq 0$, there exist a $\F_t$-measurable random variable $\tau_t$ such that $\tau_t = \tau$ on the set $\seq{\tau < t}$.
\ed

We first introduce some quantities that are specific to the study of random times. For an arbitrary random time $\tau$, we set $H:=\I_{[\![\tau,\infty)}$
and define \hfill\break 
\noindent $\bullet$ the supermartingale $Z$ associated with $\tau$, $Z :=\,^{o}(\I_{[\![0,\tau[\![})= 1- \,^{o}H$,\hfill\break 
\noindent $\bullet$ the supermartingale $\widetilde Z$ associated with $\tau$, $\wt Z :=\,^{o}(\I_{[\![0,\tau]\!]})= 1- \,^{o}(H_-)$,\hfill\break 
\noindent $\bullet$ the martingale $m :=1-\left ( \,^{o}H-H^o\right )$.

In the literature, the process $Z$ is often termed the Az\'ema supermartingale. Here we shall name the process $\wt Z$, the Az\'ema optional supermartingale, and the process $1-Z$, the Az\'ema submartingale. From the above, one can deduce that the following relationships hold:
\begin{equation}
\label{relation}
Z=m-H^o \quad \textrm{and} \quad \widetilde Z=m- (H^o)_-
\end{equation}
and we have $\wt Z - Z = \Delta H^o$, $\wt Z_+ = Z$ and $\wt Z_- = Z_-$. From Theorem 5.22 \cite{HWY}, the dual optional projection $H^o$ is of integrable variation since $H$ is of integrable variation. At time equal to zero, we have $1-\wt Z_0 =0$ and $1-\wt Z_{0+} = 1-Z_0 = 1-(\Delta H^o)_0 = 1-H^o_0$. We set $R := \inf \seq{s: Z_s = 0}$ and for a random time $\tau$, from Lemma 1.51 in \cite{AJ}, we have $\tau \leq R$.

The process $\wt Z$ is a bounded positive optional supermartingale with it's Doob-Meyer-Mertens-Gal'\v cuk decomposition given by $\wt Z = m- (H^{o})_-$. For notational simplicity and to be consistent with the notation later used in the additive decomposition of optional semimartingales, we will set $A:= (H^o)_-$. Note that $A$ is a left continuous process and to which one can apply the decomposition in \eqref{Xg} to obtain the additive decomposition $\wt Z = m - A^c - A^g$. 

\subsection{Existence of the representations}\label{s2.1}
The first main result of this paper is presented in \cref{t1}, where we obtain the existence of an additive and multiplicative representation of $\wt Z$ in terms of the drawdown and the relative drawdown of some optional supermartingale. To better understand the limitations of existing results, we give below two counter examples of finite honest times for which the Az\'ema supermartingale is not of the form given in \eqref{eqa} and \eqref{eqm}, and thus not covered by existing results in \cite{NY2, NP, K, BI, S}.

\bex\label{e01}
Let $\ff$ be the Brownian filtration and $\sigma$ a finite $\ff$-stopping time. Then $\sigma$ is an example of a finite honest time or last passage time for which represenation \eqref{eqa} and \eqref{eqm} does not hold, since otherwise it will contradict the fact that all martingales are continuous in $\ff$. 
\eex

\bex\label{exp2.2}
Let $\ff$ be the Brownian filtration and $\tau := \sup\{s: X_s = \overline{X}_s\}$ where $X_t = e^{-\sigma^2t/2 + \sigma W_t}$. The process $X$ is a positive local martingale for which $X_\infty= 0$, and by using the fact that $\sup_{t\geq 0 } \left(\sigma W_t-\sigma^2t/2\right) \sim \exp(\sigma^2)$, one can show that $\P(\tau > t\,|\,\F_t) = X_t/\overline{X}_t$ which is of the form given in \eqref{eqm}. On the other hand, it can be shown, by checking directly \cref{honest}, that for any $T\in \rr_+$ the random time $\tau' = \tau \vee T$ is a finite honest time and 
\begin{align*}
\P(\tau' > t \,|\, \F_t) = 1- \P(\tau\vee T \leq t \,|\, \F_t) = 1-\left(1-X_t/\overline{X}_t\right)\I_{\{T \leq t\}}\\
\P(\tau' \geq t \,|\, \F_t) = 1- \P(\tau\vee T < t \,|\, \F_t) = 1-\left(1-X_t/\overline{X}_t\right)\I_{\{T < t\}}
\end{align*}
which are both discontinuous at $T$. Clearly $\P(\tau' > t \,|\, \F_t)$ cannot be of the form given in \eqref{eqm}, since otherwise it will contradict the fact that all martingales are continuous in the Brownian filtration. 

\eex

\bd\label{A}
A random time $\tau$ is said to avoid all $\ff$-stopping times if for each $\ff$-stopping time $\sigma$, $\P(\tau = \sigma<\infty) = 0$.
\ed
The class of finite honest times considered in \cite{BI, K, NP, NY2, S} are assumed to avoid all $\ff$-stopping times, while in \cref{e01} and \cref{exp2.2}, the honest time $\sigma$ is a stopping time and the honest time $\tau'$ does not avoid $T \in \rr_+$. 
We mention that, as the process $\wt Z$ is not c\`adl\`ag, the standard c\`adl\`ag Skorokhod reflection lemma cannot be applied to obtain the existence and uniqueness of the additive representation as done in \cite{NY2} and the standard Doob maximal identity cannot be applied to obtain the multiplicative representation as done in \cite{BI,K, NP, NY2}. 

To find a multiplicative representation of $\wt Z$, one needs to first find the multiplicative decomposition of $\wt Z$. However, in general, the process $\wt Z$ is not c\`adl\`ag and the multiplicative decomposition of $\wt Z$ is not available in the literature. Hence we study below the multiplicative decomposition of $\wt Z$ under the assumption that $Z>0$.

\bl\label{Yfinite} Given a finite honest time $\tau$ such that $Z > 0$, then for all $t\geq 0$ the process 
\begin{align}
Y_t &:= \int_{(0,t]} \wt Z^{-1}_sdA^c_s + \int_{[0,t)} \wt Z^{-1}_{s+} dA^{g}_{s+}\label{Y}
\end{align}
is a finite l\`ag\`ad increasing process, and the optional stochastic exponential of $Y$ denoted by $\wt D = \mathcal{E}(Y)$ is an increasing process such that  both $d\wt D^c$ and $d\wt D^g_+$ are carried on the set $\{\wt Z = 1\}$.
\el 

\begin{proof}
From \eqref{Y}, we have $Y^c_t = \int_{(0,t]} \wt Z^{-1}_sdA^c_s$ and $Y^g_t =\int_{[0,t)} \wt Z^{-1}_{s+} dA^{g}_{s+}$, and the aim is to show that these integrals are finite for all $t\geq 0$. We stress that the limit of integration in the second integral is $[0,t)$ and thus the process $Y$, once shown to be finite for all $t\geq 0$, is an increasing l\`agl\`ad process. 
By (iii) of \cref{p1.1}, both the processes $Y^c$ and $Y^g$ are increasing and stopped after $\tau$. For both $Y^c$ and $Y^g$ to be well defined finite increasing processes, we need to check that the they are finite before $\tau$ and have finite right limit at $\tau$.  

We first note that, by \cref{l1.2} (i) and \cref{l1.3}, the measure $dA_+$ is carried on $\seq{\wt Z= 1} \subset \llb 0, \tau \rrb$, therefore we have $Y^c = A^c < \infty$. Secondly, for $Y^g$, we see that for $t>\tau$
\begin{gather*}
\int_{[0,t)} \wt Z^{-1}_{s+} dA^{g}_{s+} = \sum_{0\leq s\leq \tau}Z^{-1}_s\Delta^+ A^{g}_s < \infty,
\end{gather*}
which is due to the fact that for almost all $\omega$ the integrand $Z^{-1}$ is bounded away zero as $R = \inf \seq{s: Z_s = 0} = \infty$ (see Theorem 2.62 \cite{HWY}) and $A^g < \infty$.

Having shown that $Y$ is a finite increasing process, by using \cref{stochexp}, we define $\wt D$ as the optional stochastic exponential of $Y$. That is $\wt D$ is the unique solution to the following equation
\begin{align}
\wt D_t &= 1 + \int_{(0,t]}\wt D_{s-} dY^{c}_s + \int_{[0,t)}\wt D_s dY^{g}_{s+}\nonumber \\
  &= 1+ \int_{(0,t]}\wt D_{s-}\wt Z^{-1}_s dA^c_s + \int_{[0,t)}\wt D_s\wt Z^{-1}_{s+} dA^{g}_{s+} =: 1+ \wt D^c_t + \wt D^g_t \label{Dcg}
\end{align}
and it clear from \eqref{Dcg} and \cref{l1.3} that $d\wt D^c$ and $d\wt D^g$ are carried on the set $\{\wt Z = 1\}$.
\end{proof}

\bl\label{mdz}
Suppose $Z > 0$ then the process $\wt M = \wt D\wt Z$ is a c\`adl\`ag local martingale.
\begin{proof}
By an application of the It\^o formula given in \cref{ito}, 
\begin{align*}
\wt D_t\wt Z_t - \wt D_0\wt Z_0  & = \int_{(0,t]}\wt D_{s-}d\wt Z^r_s + \int_{[0,t)}\wt D_s d\wt Z^g_{s+}\\
& \quad  + \int_{(0,t]}\wt Z_{s-}d\wt D^r_s + \int_{[0,t)}\wt Z_{s} d\wt D^g_{s+} + \sum_{0\leq s<t} \Delta \wt D^g_s \Delta^+\wt Z_s + \sum_{0< s\leq t} \Delta \wt D^r_s \Delta\wt Z_s\\
   		   & = \int_{(0,t]}\wt D_{s-}d\wt Z^r_s + \int_{[0,t)}\wt D_s d\wt Z^g_{s+} + \int_{(0,t]}\wt Z_sd\wt D^c_s + \int_{[0,t)}\wt Z_{s+} d\wt D^g_{s+} = \int_{(0,t]}\wt D_{s-}dm_s 
\end{align*}
where in the last equality, we have used the fact that $\wt Z^r = m - A^c$, $\wt Z^g = -A^{g}$ and \eqref{Dcg}.
\end{proof}
\el 
From \cref{Yfinite} and \cref{mdz} we see that, under the assumption that $Z >0$, the multiplicative decomposition of $\wt Z$ is given by $\wt Z = \wt M/\wt D$. Note that, by using \cref{stochexp} we see that the unique solution to the optional stochastic exponential $\wt D$ is given by 
\begin{align*}
\wt D = e^{Y^c}e^{Y^g}\prod_{0\leq s< \cdot} (1+\Delta^+ Y_s)e^{\Delta^+ Y_s}.
\end{align*}
The process $\wt D$ can then be further decomposed multiplicatively into it's continuous and left continuous parts. That is $\wt D = D^cD^g$, where
\begin{align}
D^c &:= e^{Y^c} \quad \mathrm{and}\quad D^g := e^{Y^g}\prod_{0\leq s< \cdot} (1+\Delta^+ Y_s)e^{\Delta^+ Y_s}.\label{D}
\end{align}
Here we recall that 
\bde
Y^c_t  = \int_{(0,t]} \wt Z^{-1}_sdA^c_s = A^c   \quad \mathrm{and} \quad Y^g_t =\int_{[0,t)} \wt Z^{-1}_{s+} dA^{g}_{s+}.
\ede 
From the form of $Y^c$ and $Y^g$, we see that $D^c$ and $D^g$ are strictly positive increasing processes such that $dD^c$ and $dD^g_+$ are carried on the set $\seq{\wt Z = 1}$.

We are now in a position to study the additive and multiplicative representations. The key idea is that, instead of the local martingales $m$ and $\wt M$, we consider the local optional supermartingales 
\begin{align}\label{nN}
n : = m - A^{g} = \wt Z + A^c \qquad \mathrm{and} \qquad N: = \wt ZD^c = \wt Z e^{A^c}.
\end{align}
We stress here that the process $N: = \wt ZD^c$ is well defined, even if $Z$ is not strictly positive, since $D^c = e^{A^c}$ is always well defined. 
In general, the processes $n$ and $N$ are not necessarily c\`adl\`ag on $\llb 0 ,\tau \rrb$ and due to \cref{p1.1} they are only c\`adl\`ag on $\rrb \tau ,\infty \llb$. We remark that since $A^{g}$ has only positive jumps, the processes $\overline{n}$ and $\overline{N}$ must be c\`adl\`ag and hence optional processes. Therefore the sets $\seq{n= \overline{n}}$ and $\seq{N = \overline{N}}$ are optional sets.

\bl\label{l2.1}
For any finite honest time $\tau$, we have $\seq{n = \overline{n}} = \seq{\wt Z = 1} = \seq{N =  \overline{N}}$.
\begin{proof}
From \cref{p1.1}, we observe that
\be\label{teq1}
\tau = \sup\seq{s: \wt Z_s = 1} = \sup\seq{s:  n_s=1+ A^c_{s} }.
\ee 
From the inequality $\wt Z \leq 1$, we deduce that $n \leq 1+ A^c$ and $\seq{n=1+ A^c}  \subseteq  \seq{n=\overline n}$. Using \eqref{teq1} and the fact that the process $1+A$ is constant after $\tau$ (from \cref{l1.2} (i) and \cref{l1.3}), the process $\overline{n}$ must equal to the constant process $1+A^c$ after $\tau$. This together with the fact that $\seq{\wt Z = 1}$ is contained on $\llb 0, \tau \rrb$, we have
\bde
\seq{n=\overline{n}}\,\cap \,\rrb \tau ,\infty \llb = \seq{\wt Z=1}\,\cap \,\rrb \tau ,\infty \llb = \emptyset.
\ede 
This implies $\seq{n=1+ A^c}  \subseteq  \seq{n=\overline n} \subseteq \llb 0,\tau\rrb$. By \cref{l1.2}, the set $\seq{\wt Z = 1}$ is the largest optional set contained in $\llb 0, \tau\rrb$ from which we conclude that $\seq{n= \overline {n}} = \seq{\wt Z= 1}$. Similar arguments shows that $\seq{N= \overline {N}} = \seq{\wt Z= 1}$
\end{proof}
\el 

\brem 
The set equality $\seq{n = \overline{n}} = \seq{n = 1+A^c}$ implies that $\overline{n} = 1+A^c$ on the set $\seq{\wt Z = 1}$. Intuitively, the equality $\overline{n} = 1+A^c$ should also hold everywhere, since they have the same initial condition and both $d\overline{n}$ and $dA^c$ are carried on the set $\seq{\wt Z = 1}$. 
\erem 
To prove the above observation and therefore the additive and multiplicative representation of $\wt Z$, we make use of the following time change process,
\begin{align*}
g_t   &= \sup\seq{s \leq   t: \wt Z_s = 1}.
\end{align*}
Note that $\wt Z_\tau = 1$, but in general, it is not true that $\wt Z_{g_t} = 1$.

\bt\label{t1}
Let $\tau$ be a finite honest time random time.\\
\noindent (i) An additive representation of $\wt Z$ is given by
\begin{align*}
\wt Z &= 1+ n - \overline{n},
\end{align*}
where $n = \wt Z + A^{c}$ and $1+A^c = \overline {n}$.\\
\noindent (ii) A multiplicative representation of $\wt Z$ is given by
\begin{align*}
\wt Z & = N/\overline{N}
\end{align*}
where $N = \wt ZD^c$ and $D^c =e^{A^c}= \overline {N}$.
\begin{proof}
The goal of the proof is to show that $\overline{n} = 1+A^c$ and $\overline{N} = e^{A^c}$. Firstly, the processes $\overline{n}$ and $1+A^c$ have the same initial condition and it is clear that $\overline{n} \leq 1+A^c$. To show the reverse inequality, we must consider two cases. Given a finite stopping time $T$, we first suppose that $(\omega,g_T(\omega)) \in \seq{\wt Z=1}$, then by using the facts that $A^c$ is continuous, $dA^c$ is carried on the set $\{\wt Z = 1\}$ and, by Lemma \ref{l2.1}, that $\seq{n = \overline{n}} = \seq{n = 1+A^c} = \seq{\wt Z = 1}$, we have $1+ A^c_T  = 1+ A^c_{g_T} = \overline{n}_{g_T} \leq \overline{n}_T$. On the other hand, suppose $(\omega,g_T(\omega))\not \in \seq{\wt Z = 1}$ but belongs to the right closure of $\seq{\wt Z = 1}$, then there exists an increasing sequence of random times $(g^n_T)_{n\in\nn}$ in $\seq{\wt Z = 1}$ such that $g^n_T\uparrow g_T$ a.s.. Then similar to the previous case, we have $1+ A^c_T  = 1+ A^c_{g_T} = \overline{n}_{g_T-} \leq \overline{n}_T$, which implies $1+ A^c = \overline{n}$.

For the multiplicative representation, we first observe that $dD^{c}$ is carried on the set $\seq{\wt Z= 1} = \seq{N= D^c}$. Then one can repeat the arguments used in the proof of the additive representation with $N$, $\overline{N}$ and $D^c$ in place of $n$, $\overline{n}$ and $1+ A^c$ to conclude that $\overline{N}$ is equal to $D^c$. 
\end{proof}
\et

\bex\label{ex2.1}
To illustrate our main result obtained in \cref{t1}, let us re-visit below the honest times given in \cref{e01} and \cref{exp2.2}. We recall from \cref{exp2.2} that the Az\'ema optional supermartingale $\wt Z$ associated with the finite honest time $\tau'$ is 
\begin{gather*}
\wt Z_t = \P(\tau' \geq t \,|\, \F_t) = 1-\left(1-X_t/\overline{X}_t \right)\I_{\{T < t\}}.
\end{gather*}
By applying the It\^o's formula to $\wt Z$, using the uniqueness of the Doob-Meyer-Mertens-Gal'\v cuk decomposition, and the fact that $d\overline{X}$ is carried on the set $\{X = \overline{X}\}$, one can deduce that 
\begin{align*}
A^c_{t} = \ln(\overline X_{t\vee T}) - \ln(\overline X_{T}).
\end{align*}
From \eqref{D} and \eqref{nN} we have $D^c_t = e^{A^c_t} = \overline{X}_{t\vee T}/\overline{X}_{T}$ and $N_t = \wt Z_t D^c_t = \left(\I_{\{T \geq t\}} + (X_t/\overline{X}_T)\I_{\{T < t\}}\right)$. From the form of $N$ we clearly have $\overline{N}_t = \overline{X}_{t\vee T}/\overline{X}_{T} = D^c_t$ and therefore $\wt Z_t = N_t/\overline {N}_t$ for all $t\geq 0$.

On the other hand, the multiplicative decomposition for the Az\'ema optional supermartingale associated with a stopping time $\sigma$, given in \cref{e01}, is trivial in the sense that $\wt Z = \I_{\llb 0, \sigma \rrb} = 1 - A^g$ and clearly $\wt Z = N/\overline{N}$ where $N = \wt Z =\I_{\llb 0, \sigma \rrb}$ is an optional supermartingale.

\eex

\brem 
To obtain non-trivial examples, such as $\tau'$ studied in \cref{exp2.2}, we need to find finite honest times for which both $A^c$ and $A^g$ are non-zero. Honest times with such property can be constructed by taking known examples of finite honest times $\tau$ which avoids all stopping times and consider the honest time $\tau \vee \sigma$ where $\sigma$ is any finite stopping time. We will discuss this type of construction in more details in \cref{e3.1} once we have developed some generic tools in \cref{os}.

\erem

\subsection{Characterisation of honest times \& uniqueness of the representations}\label{sec2.1}

The aim of this section is of twofold. We first extend the {\it Doob maximal identity} in \cref{dme} in order to provide in \cref{cor2.1} a characterisation of finite honest times using optional supermartingales of class $\mathcal{N}_0$, defined in \cref{n0}. Then we prove the uniqueness of the multiplicative representation in \cref{p2.18} and, by using a l\`agl\`ad variant of the {\it Skorokhod reflection lemma} obtained in \cref{srl}, we prove the uniqueness of the additive representation in \cref{p2.21}.

\bd\label{n0}
An local optional supermartingale $N$ is said to belong to the class $\mathcal{N}_0$ if \\
(i) the process $N$ is non-negative and $\lim_{t\rightarrow \infty} N_t = 0$, \\
(ii) the running supremum $\overline{N}$ is continuous,\\
(iii) the graph of $\tau := \sup\seq{s: N_s = \overline{N}_s}$ belongs to $\{N = \overline{N}\}$ or equivalently $N_{\tau} = \overline{N}_{\tau}$,\\
(iv) the process $N$ exhibits the decomposition $N = N_0 + M^N - A^N$ where $M^N$ is a local martingale and $A^N$ is a left continuous increasing process such that $dA^N_+$ is carried on $\seq{N = \overline{N}}$.
\ed 
The class $\mathcal{N}_0$ extends the notion of local martingales of class $\mathcal{C}_0$ and class $\mathcal{M}_0$, where the class $\mathcal{M}_0$ consists of non-negative local martingales with continuous running supremum which converges to zero at infinity, see page 616 of \cite{NP}, and $\mathcal{C}_0$ consists of continuous local martingales of class $\mathcal{M}_0$.

\bd
An local optional supermartingale $N$ is said to belong to the class $\mathcal{N}^*_0$ if $N \in \mathcal{N}_0$, $N_0 = 1$ and $A^N$ is a pure jump process.
\ed 

\brem 
The process $N :=\wt ZD^c$ given in \cref{t1} belong to class $\mathcal{N}_0^*$. The process $N = \wt ZD^c$ is clearly non-negative and $\lim_{t\rightarrow \infty} N_t = 0$ since $\tau$ is finite and hence $\lim_{t\rightarrow \infty} \wt Z_t = 0$. We know from \cref{p1.1} that $\wt Z_\tau = N_\tau/\overline{N}_\tau = 1$ and by applying the It\^o formula in \cref{ito} to $N = \wt ZD^c$, we have
\begin{align*}
D_t^c\wt Z_t 
		   & = 1+ \int_{(0,t]}D^c_{s-}d\wt Z^r_s + \int_{[0,t)} D_s^c d\wt Z^g_{s+} + \int_{(0,t]}\wt Z_sd\wt D^c_s = 1+ \int_{(0,t]} e^{A^c_s} dm_s - \int_{[0,t)}e^{A^c_s} dA^g_{s+}. 
\end{align*}
From the above, we see that $A^N$ is a left continuous increasing pure jump process and $dA^N_+$ is carried on $\{N = \overline{N}\}$ since $dA^g_+$ is carried on $\{\wt Z=  1\}$ which is equal to $\{N = \overline{N}\}$ by \cref{l2.1}.
\erem

\bl [Variant of the Doob maximal identity]\label{dme}
Suppose $N$ is an local optional supermartingale which belongs to $\mathcal{N}_0$ and we let $\tau := \sup\seq{s: N_s = \overline{N}_s}$. Then $\tau$ is a finite honest time such that its Az\'ema optional supermartingale is given by $\wt Z = N/\overline{N}$.
\el

\begin{proof}
It is clear that $\tau$ is a finite last passage time and hence a finite honest time. Therefore we need only to compute $\wt Z$ associated with $\tau$. Let us consider the process $Y = 1- N/\overline{N}$ which is positive and bounded by one. Hence $Y$ is a positive optional submartingale of class-$(D)$ and $Y_\infty = 1$. By applying the It\^o formula to $Y$ and using the uniqueness of the Doob-Meyer-Mertens-Gal'c\v uk decomposition of $Y$ given in \cref{dm1}, we can conclude that
\begin{gather*}
M^Y_t =  -\int_{(0,t]}\overline{N}^{-1}_s dM^N_s  \quad \mathrm{and} \quad A^Y_t = \int_{[0,t)}\overline{N}^{-1}_s dA^N_s + \ln(\overline{N}_t).
\end{gather*}
where $M^Y$ optional martingale and hence a uniformly integrable martingale, and $dA^Y_+$ is carried on $\{Y = 0 \}= \{N = \overline{N} \}$. Let $\gamma_t := \inf\seq{s\geq t: N = \overline{N}} = \inf\seq{s\geq t: Y = 0}$, which by convention takes the value infinite if the set is empty. We observe that for every stopping time $T$,
\begin{align}
Y_{\gamma_T} = Y_{\gamma_T}\I_\seq{\tau < T} + Y_{\gamma_T}\I_\seq{\tau \geq T} = \I_\seq{\tau < T}. \label{xgammaT}
\end{align}
To obtain the second equality above, we notice that $Y_{\gamma_T}\I_\seq{\tau < T} = Y_\infty\I_\seq{\tau < T} = \I_\seq{\tau < T} $ and 
\begin{gather*}
Y_{\gamma_T}\I_\seq{\tau \geq T} =  Y_{\gamma_T}\I_\seq{\tau > T} + Y_{\gamma_T}\I_\seq{\tau =  T}.
\end{gather*}
On the set $\{\tau > T\}$, the equality $N_{\gamma_T} = \overline N_{\gamma_T}$ clearly holds for points $(\omega,\gamma_T(\omega)) \in \seq{N = \overline{N}}$. Suppose now that $(\omega,\gamma_T(\omega))\not \in \seq{N = \overline{N}}$ but is in the left closure of the set $\{N = \overline{N}\}$. That is there exists a decreasing sequence of stopping times $(\gamma^n_T)_{n\in \nn}$ such that $\gamma^n_T \downarrow \gamma_T$ almost surely and for every $n$, $(\omega, \gamma^n_T(\omega)) \in \seq{N = \overline{N}}$. This together with the continuity of $\overline{N}$ implies that $N_{\gamma_T+} = \overline{N}_{\gamma_T} > N_{\gamma_T}$. However, since $\Delta^+ N = -\Delta^+ A^N\leq 0$, we must have $N_{\gamma_T} \geq N_{\gamma_T+}$ which is a contradiction. While, on the set $\{\tau =  T\}$, we have 
\begin{gather*}
N_{\gamma_T} = N_{\tau}\I_{\{N_\tau = \overline N_\tau\}} + N_\infty\I_{\{N_\tau < \overline N_\tau\}} = \overline{N}_{T} = \overline{N}_{\gamma_T},
\end{gather*}
where we have used the fact that $N_\infty = 0$, $\gamma_T = \gamma_\tau = \tau = T$ on the set $\{N_\tau = \overline N_\tau\}$ which is assumed to be of probability one, and the fact that $\overline{N}$ is continuous. This shows that $Y_{\gamma_T}\I_\seq{\tau \geq T}  = 0$.

By taking the $\F_T$ conditional expectation of both hand sides of \eqref{xgammaT}, we obtain
\begin{align*}
\mathbb{E}(Y_{\gamma_T}\,|\,\F_T) = 1-\wt Z_T.
\end{align*}
Using the fact that $M^Y$ is uniformly integrable, we have from the Doob optional sampling theorem, $\mathbb{E}(Y_{\gamma_T}\,|\,\F_T) = M^Y_T + \mathbb{E}(A^Y_{\gamma_T}\,|\,\F_T)$. Finally, as $A^Y$ is left-continuous and $dA^Y_+$ is carried on $\{N = \overline{N}\}$, we have $A^Y_{\gamma_T} = A^Y_{T}$ and hence $\wt Z_T = N_T/\overline{N}_T$ for all finite stopping times $T$.
\end{proof}

\bcor\label{cor2.1}
Suppose that $\tau$ is a finite honest time then there exists an optional supermartingale $N$ of class $\mathcal{N}_0$ such that $\tau$ is the end of the optional set $\seq{N = \overline{N}}$ and $\wt Z = N/\overline{N}$. Conversely, given a local optional supermartingale $N$ of class $\mathcal{N}_0$, the end of the optional set $\seq{N = \overline{N}}$ is a finite honest time such that $\wt Z = N/\overline{N}$.

\begin{proof}
It is sufficient to combine \cref{l2.1}, \cref{t1} and \cref{dme}.
\end{proof}
\ecor

The above corollary gives a characterization of the finite honest times through local optional supermartingales.
However, given a finite honest time $\tau$, the class $\mathcal{N}_0$ is too big to have the uniqueness of the multiplicative representation. It is not hard to see that $\wt Z = N/\overline{N}$ where $N$ can be either $\wt Z$, $\wt ZD^c$, $k\wt Z$ or $k\wt ZD^c$ for any $k\geq 0$, which all belongs to the class $\mathcal{N}_0$. The class $\mathcal{N}^*_0$ is introduced to restrict our attention to the case $k = 1$ and to remove, whenever possible, the trivial candidate $\wt Z$.

\bp\label{p2.18}
Suppose that $\tau$ is a finite honest time then there exists a unique optional supermartingale $N$ of class $\mathcal{N}^*_0$ such that $\tau$ is the end of the optional set $\seq{N = \overline{N}}$ and $\wt Z = N/\overline{N}$.
\ep

\begin{proof}
We need only to show the uniqueness of the process $N$ inside the class $\mathcal{N}^*_0$. Suppose that there exist another process $X \in \mathcal{N}^*_0$ such that $\wt Z = N/\overline{N} = X/\overline{X}$ and the decomposition of $N$ and $X$ are given by $N = N_0 + M^N - A^N$ and $X = X_0 + M^X - A^X$ where $A^N$ and $A^X$ are left continuous increasing pure jump processes. Then by the l\`agl\`ad It\^o formula and the fact that $d\overline{N}$ and $d\overline{X}$ are carried on the set $\{N = \overline{N}\}$ and $\{X = \overline{X}\}$ respectively, we have
\begin{align*}
\wt Z_t & = 1+ \int_{(0,t]}\overline{X}^{-1}_s dM^X_s - \int_{[0,t)}\overline{X}^{-1}_s dA^X_s - \ln(\overline{X}_t) \\
		& = 1+ \int_{(0,t]}\overline{N}^{-1}_s dM^N_s - \int_{[0,t)}\overline{N}^{-1}_s dA^N_s - \ln(\overline{N}_t).
\end{align*}
From the uniqueness of Doob-Meyer-Merten-Gal'\v cuk decomposition of $\wt Z = n - A$. We deduce that 
\begin{align*}
n & = 1+ \int_{(0,t]}\overline{X}^{-1}_s dM^X_s = 1+ \int_{(0,t]}\overline{N}^{-1}_s dM^N_s\\
A_t & = \int_{[0,t)}\overline{X}^{-1}_s dA^X_{s+} + \ln(\overline{X}_t) = \int_{[0,t)}\overline{N}^{-1}_s dA^N_{s+} + \ln(\overline{N}_t).
\end{align*}
From the continuity of $\overline{X}$ and $\overline{N}$ and the fact that $A^X$ and $A^N$ are left continuous pure jump processes, we deduce that $A^c_t = \ln(\overline{N}_t) = \ln(\overline{X}_t)$ which implies that $\overline{N}_t = \overline{X}_t$ and $A^g_t  = \int_{[0,t)}\overline{N}^{-1}_s dA^N_{s+} = \int_{[0,t)}\overline{X}^{-1}_s dA^X_{s+}$ and consequently $A^N = A^X$. Finally, since $\overline{N}_t = \overline{X}_t$, we can conclude that $M^X= M^N$ and hence $X - X_0 = N - N_0$, where $X_0 = N_0 = 1$.
\end{proof}

\brem
The class of finite honest times for which the process $N \in \mathcal{N}_0^*$ is the trivial candidate $\wt Z$ is exactly the class of finite {\it thin honest times}. That is honest times whose graph is contained in the disjoint union of the graph of a family of $\ff$-stopping times or equivalently the process $H^o_- = A$ is a pure jump increasing process (see Definition 1.1 and Theorem 1.4 in \cite{ACJ}). In this case, we have $A = A^g$, $A^c = 0$, $D^c = 1$ and $N= \wt ZD^c = \wt Z$.

The counter example mentioned in the introduction, given by Acciaio and Penner \cite{BI}, is an example of a thin honest time, and the multiplicative representation holds trivially with $N = \wt Z$. As this example is quit involved, we refer to readers to Proposition 4.8 of \cite{ACDJ}. For a simpler illustration, let us return to \cref{e01} and consider a finite $\ff$-stopping time $\sigma$, which is a finite thin honest time with $\wt Z = \I_{\llb 0, \sigma\rrb}$. Suppose that there exist $N \in \mathcal{N}_0$ such that $\wt Z = N/\overline{N}$, then we can deduce from the equality $\overline {N} \I_{\llb 0, \sigma\rrb} =  N$ that $N$ must be non-decreasing on $\llb 0, \sigma\rrb$ and zero on $\rrb \sigma, \infty \rrb$. In fact, the process $N$ must be $k\I_{\llb 0, \sigma\rrb}$ from some $k\geq 0$ since $N$ is an local optional supermartingale. If we restrict ourselves to the class $\mathcal{N}_0^*$ then it is necessary that $k = 1$ and $N = \wt Z$.
\erem 

We now investigate the uniqueness of the additive representation of $\wt Z$. To do this, we provide a l\`agl\`ad variant of the Skorokhod reflection lemma which is not available in the literature. 
\bl[Variant of the Skorokhod reflection lemma]\label{srl}
Let $y$ be a real-valued l\`agl\`ad function on $[0,\infty)$ such that $y(0) = 0$ and it's running infimum $\inf y = \underline{y}$ is continuous. Then, there exists a unique pair $(z, a)$ on $[0,\infty)$ where $a(t)=\sup_{s \le t}-y(s)$, satisfying
the following conditions:\\
(i) $z(t) = y(t) + a(t) \geq 0$ for all $t \geq 0$,\\
(ii) $a$ is a increasing, continuous function with initial value zero,\\
(iii) the measure da is carried on the set $\{t : z(t) = 0\}$.
\el
\begin{proof}
See \cref{4.3}.
\end{proof}

\bp\label{p2.21}
Given a finite honest time $\tau$, there exists a unique local optional supermartingale $n$ with $n_0 = 1$ and continuous running supremum, such that $\wt Z = 1 + n- \overline n$.
\ep

\begin{proof}
The existence of $n$ with $n_0 = 1$ and continuous running supremum follows from \cref{t1}. The uniqueness of the additive representation of $\wt Z$ follows from \cref{srl} if we set $y= 1- n$, which gives $a = \overline{-1+n} = -1+ \overline{n}$ and $z = 1-n + a  =  -n+ \overline{n}$. 
\end{proof}

\brem\label{r2.6}
We stress that, in Theorem \ref{t1}, the core of the proof is showing that $\overline{N} = e^{A^c}$ and $\overline{n} = 1+A^c$, which gives the important property that both $\overline{N}$ and $\overline{n}$ are continuous.  In fact, one can argue that the continuity of both $\overline{N}$ and $\overline{n}$, where $n$ and $N$ are defined in \eqref{nN}, is the most important property. Since if one can show directly that $\overline{N}$ and $\overline{n}$ are continuous then the existence and uniqueness of the additive representation of $\wt Z$ can be obtained through \cref{srl}. While the the existence and uniqueness of the multiplicative representation of $\wt Z$ can be obtained by combining \cref{l2.1}, \cref{dme} and \cref{p2.18}.
\erem
We conclude the first part of the paper by showing, for the sack of completeness, that the processes $\overline{n}$ and $\overline{N}$ are continuous without showing that they are equal to $1+A^c$ and $e^{A^c}$. This gives us an alternative method to prove the existence of the additive and the multiplicative representation of $\wt Z$ and highlights the importance of \cref{l2.1}.
\bp
Given a finite honest time $\tau$, the running supremum of the processes $n = m-A^g$ and $N =  \wt ZD^c$ are continuous.
\begin{proof}
We will only present the proof of continuity for $\overline{n}$ since the proof of continuity for $\overline{N}$ follows from similar arguments. To this end, suppose that the supremum of $n$ is not continuous and the left jumps of $n$, that is the jump of the martingale $m$, can take $n$ to its supremum. More specifically, we set $T:= \inf\seq{s: \Delta \overline{n}_s > 0}$ and suppose that $T<\infty$. Using the fact that $\overline n$ is c\`adl\`ag, we deduce that $T> 0$ and $\llb T \rrb \in \seq{n = \overline{n}}$ which by \cref{l2.1} is equal to $\seq{\wt Z = 1}$. Then for fixed $\omega \in \Omega$, there are two cases to consider, (i) the point $T(\omega)$ is a left isolated point of the set $\seq{s:\wt Z_s(\omega) =1}$ and (ii) the point $T(\omega)$ is not a left isolated point of the set $\seq{s:\wt Z_s(\omega) =1}$. 

In case (i), we consider the random time $\tau_t = \sup\seq{s <   t: \wt Z_s = 1}$. Note that since $T(\omega)$ is a left isolated point and $dA^c$ is carried on the set $\seq{n = 1+A^c} = \seq{\wt Z = 1} = \seq{n = \overline{n}}$, we must have, for the fixed $\omega$, $\tau_T < T$ and $A^c_{\tau_T} = A^c_{T-} =A^c_{T}$. However this is a contraction since this will imply
$$1+A^c_{T-} = 1+A^c_{\tau_T}  = \overline{n}_{\tau_T} = \overline{n}_{T-} < \overline{n}_{T} = 1+A^c_T.	$$
We point out that, in this case, one does not have to distinguish whether $\tau_t(\omega)$ belongs to the set $\seq{s:\wt Z_s(\omega) =1}$ or is in its right closure, since $\overline{n}$ is continuous before $T$.

In case (ii), since $T(\omega)$ is not a left isolated point of $\seq{s:\wt Z_s(\omega) =1}$, there exists an increasing sequence $(T_n(\omega))_{n \in \nn}$ such that $\forall n\in \nn$, $T_n(\omega) < T(\omega)$, $T_n(\omega) \in \seq{s:\wt Z_s(\omega) =1}$ and $T_n(\omega)\uparrow T(\omega)$. This implies that for the fixed $\omega$, $1+A^c_{T-}= \overline{n}_{T-}< \overline{n}_T = 1+A^c_{T}$ and this gives a contradiction.
\end{proof}
\ep

\section{Optional Semimartingales of Class-$(\Sigma)$}\label{os}
In this part of the paper, we study the Az\'ema supermartingale of finite honest times in a general context and extend, in \cref{d3.2}, the notion of semimartingales of class-$(\Sigma)$ to {\it optional semimartingales of class-$(\Sigma)$} by allowing for jumps in the finite variational part of the semimartingale decomposition. The goal below is to recover some existing results in the literature for semimartingales of class-$(\Sigma)$ in the context of optional semimartingales of class-$(\Sigma)$, and to apply them to the construction of finite honest times. Although some results presented below might not be surprising, we believe that the techniques used are of interest as we no longer deal with c\`adl\`ag processes. 

To be specific, we extend Lemma 2.2 (1)-(3), Lemma 2.3 and Lemma 2.4 from \cite{CNP} in \cref{l3.1}, \cref{l3.3} and \cref{l3.2} respectively. Secondly, we extend Theorem 3.1 (1) in \cite{CNP} by showing in \cref{t3.1} and \cref{c3.1} that {\it Madan-Roynette-Yor} type formulae, which relates the price of a put/call option with the last passage time of zero of the pay-off, can be recovered for optional semimartingales of class-$(\Sigma)$. Lastly, by using \cref{l3.3} and \cref{c3.1}, we obtain in \cref{max} a method to construct finite honest times for which the multiplicative decomposition of $\wt Z$ obtain in \cref{t1} is non-trivial, in that $N \neq \wt Z$.

\bd\label{d3.1}
An optional semimartingale $X$ with decomposition $X = X_0+ M + A$ where $M$ is a local martingale with $M_0= 0$ and $A = A^c + A^d + A^g$ is a l\`agl\`ad process of finite variation with $A_0 =0$ is said to satisfy the Skorokhod minimal reflection condition at zero if for every $t\geq 0$,
\begin{gather*}
\int_{[0,t)} \I_\seq{X_s \neq 0 } (dA^c_s+dA^g_{s+}) = 0 \quad \mathrm{and} \quad \int_{(0,t]}\I_\seq{X_{s-} \neq 0 } dA^d_s= 0.
\end{gather*}
\ed

\bd\label{d3.2}
An optional semimartingale $X$ is said to be of class-$(\Sigma)$ if it satisfies the Skorokhod minimial reflection condition and $X_0 = 0$ and $A^d = 0$.
\ed

In the existing definition of semimartingale of class-($\Sigma$) given in \cite{BI,EOM, EOMT, K,N1, NY2, Y1,Y2}, the process of finite variation $A$ in the decomposition of $X$ is continuous by definition, that is $A = A^c$. 
The current extension is non-trivial in that recent studies of honest times in the Poisson filtration have provided explicit examples of positive optional submartingales of class-$(\Sigma)$ whose finite variation part $A$ is a pure jump process, i.e. $A^c= 0$ and $A^g \neq 0$. In fact, it is proven in Theorem 3.6 of Aksamit et al. \cite{ACJ}, that in any jumping filtration, for example the Poisson filtration, the finite variation part $A$ in the Doob-Meyer-Mertens-Gal'\v cuk decomposition of the Az\'ema optional supermartingale $\wt Z$ associated with a finite honest time must be a pure jump process.

\bex
Given a finite honest times $\tau$, we see from \cref{l1.3} that the process $1-\wt Z$ is a positive optional submartingale of class-($\Sigma$) such that $\llb \tau \rrb \in \{1- \wt Z = 0\}$.
\eex

From this point onwards, given an optional semimartingale $X = M + A$, the process $M$ and $A$ will denote the local martingale and the l\`agl\`ad process of finite variation in the optional semimartingale decomposition of $X$. The left jumps are given by $\Delta X = \Delta M$, the right jumps are given by $\Delta^+ X = \Delta A_+ = \Delta A^g_+$ and $\seq{\Delta^+ X \neq 0} \subset \seq{X= 0}$.

\bl \label{l3.1}
Let $X$ be an optional semimartingale of class-$(\Sigma)$ then:\hfill\break 
(i) The processes $X^+$, $X^-$ and $|X|$ are local optional submartingales.\hfill\break
(ii) If $\Delta X \geq 0$ then $X^+$ is of class-$(\Sigma)$.\hfill\break 
(iii) If $\Delta X \leq 0$ then $X^-$ is of class-$(\Sigma)$.\hfill\break 
(iv) If $\Delta X = 0$ then $|X|$ is of class-$(\Sigma)$.\hfill\break 
(v) If $X$ is a positive optional submartingale then $A^c =(\overline{-n})\vee 0$ where $n := M + A^g$.
\begin{proof}
(i) The fact that $X^+$, $X^-$ and $|X|$ are local submartingales follows directly from the Tanaka formula in \cref{tanaka2} and the fact that $dA_+$ is carried on $\seq{X= 0}$.

\noindent (ii) We prove only $(ii)$ as the proof of $(iii)$ and $(iv)$ are similar. By \cref{tanaka2} 
\begin{align*}
X_t^+ &=   
\int_{(0,t]} \I_\seq{X_{s-} > 0}d(A^c_s + M_s) + \int_{[0,t)} \I_\seq{X_{s} > 0}\, dA^g_{s+} + \sum_{0< s \leq  t}  \I_\seq{X_{s-} > 0 }(X_{s})^- \\
 & \quad + \sum_{0< s \leq  t} \I_\seq{X_{s-} \leq 0 }(X_{s})^+ + \sum_{0\leq s < t}  \I_\seq{X_s > 0 }(X_{s+})^- + \sum_{0\leq s < t}   \I_\seq{X_s \leq 0 }(X_{s+})^+  + \frac{1}{2}L^0_t(X).
\end{align*}
The process $L^0_t(X)$ is called the local time of $X$ at zero and by \cref{tanaka3}, the measure $dL^0(X)$ is carried on the set $\seq{X= 0} \subseteq \seq{X^+ = 0}$. Then by using the fact that $\Delta X\geq 0$, we have 
\begin{align*}
X_t^+  &=  \int_{(0,t]} \I_\seq{X_{s-} > 0}dM_s + \sum_{0< s \leq  t}  \I_\seq{X_{s-} \leq 0 }(X_{s})^+ + \sum_{0\leq s < t}  \I_\seq{X_s > 0 }(X_{s+})^- \\
	& \quad + \sum_{0\leq s < t}  \I_\seq{X_s \leq 0 }(X_{s+})^+  + \frac{1}{2}L^{0}_t(X).
\end{align*}
Note that the right hand jumps
$\sum_{0\leq s < t}  \I_\seq{X_s > 0 }(X_{s+})^- + \I_\seq{X_s \leq 0 }(X_{s+})^+$ are supported on the set $\seq{X^+ = 0}$ because $\seq{\mathrm{sign}(X) \neq \mathrm{sign}(X_+)} \subseteq \seq{\Delta^+X > 0} \subseteq \seq{X = 0} \subseteq \seq{X^+ = 0}$.

The l\`agl\`ad process of finite variation $A$ in the optional semimartingale decomposition of $X$ is left continuous and therefore predictable. This implies that there exists a localising sequence of stopping times $(T_n)_n$ such that $X^+ = (M + A)^+$ is integrable and $M^{T_n}$ is a uniformly integrable martingale. This implies that the increasing process
\begin{align*}
V_t = \sum_{0< s \leq  t}   \I_\seq{X_{s-} \leq 0 }(X_{s})^+
\end{align*}
stopped at $T_n$ is of integrable variation and the dual predictable projection $V^p$ of $V$ exists and is locally of integrable variation. 

To show that $V^p$ is continuous, following similar argument to Lemma 2.2 in \cite{CNP}, we note that on the set $\seq{\Delta V>0}$, the jump $\Delta V$ is bounded by $\Delta X = \Delta M\geq 0 $. Therefore 
\begin{align*}
\Delta V^p = \,^p(\Delta V) & \leq \,^p(\Delta M)
\end{align*}
and from the predictable sampling theorem, $\,^p(\Delta M)_T = 0$ for all predictable stopping times $T$.  To this end, by using the continuity of $V^p$ we obtain
\begin{align*}
\mathbb{E}\big(\int_{[0,T_n)} \I_\seq{X^+_s > 0} dV^p_s\,\big)  & = \mathbb{E}\big(\int_{[0,{T_n})} \I_\seq{X^+_{s-} > 0}  dV^p_s\,\big)\\
									    & = \mathbb{E}\big(\int_{[0,{T_n})} \I_\seq{X^+_{s-} > 0}\I_\seq{X_{s-}\leq 0} dV_s \,\big) = 0.
\end{align*}
Finally, by monotone convergence theorem, we let $n\rightarrow \infty$ to show that $V^p$ is supported on $\seq{X^+=0}$.  By similarly arguments we can conclude that $X^-$ and $|X|$ are of class-$(\Sigma)$.

(v) From the fact that $X\geq 0$, we have $-n\leq A^c$ and hence $\seq{X = 0}   = \seq{-n = A^c} \subseteq \seq{-n = \overline{-n}}$.  This shows that $A^c = \overline{-n}$ on the set $\seq{X = 0}$. It is also evident that the processes $(\overline{-n})\vee 0$ and $A^c$ have the same initial condition.  Using the the inequality $-n \leq A^c$, we can conclude that $\overline{-n} \leq A^c$ and $(\overline{-n})\vee 0  \leq A^c$. To show the reverse inequality, let $g_t = \sup\{s\leq t: X_s = 0\}$ and given any stopping time $T$, either $g_T = 0$ or $g_T>0$. In the case where $g_T(\omega) > 0$ and $(\omega,g_t(\omega)) \in \seq{X=0}$, we have from the continuity of $A^c$ that
\begin{align*}
A^c_T = A^c_{g_T} = (\overline{-n})_{g_T} = (\overline{-n})_{g_T} \vee 0 \leq (\overline{-n})_T \vee 0.
\end{align*}
In the case where $g_T(\omega) > 0$ and $(\omega,g_T(\omega))$ does not belong to ${X=0}$ but is in the right closure of $\seq{X=0}$, we have by continuity of $A^c$
\begin{align*}
A^c_T = A^c_{g_T} = (\overline{-n})_{g_T-} \leq (\overline{-n})_{g_T} \vee 0 \leq (\overline{-n})_T \vee 0.
\end{align*}
In the case where $g_T(\omega) = 0$, we have $A^c_T = A^c_0 = 0 \leq (\overline{-n})_T \vee 0$. Hence we have $A^c = (\overline{-n})\vee 0$.
\end{proof}

\el

\bl \label{l3.3}
Given two optional semimartingales $X$ and $Y$ of class-$(\Sigma)$, suppose $[M^X, M^Y]= 0$ then $XY$ is an optional semimartingale of class-$(\Sigma)$.
\el 
\begin{proof}
By application of the It\^o formula in \cref{ito}
\begin{align*}
& X_tY_t - X_0Y_0 \\
	& = \int_{(0,t]}Y_{s-}dX^r_s + \int_{[0,t)}Y_s d X^g_{s+} + \int_{(0,t]} X_{s-}dY^r_s + \int_{[0,t)} X_{s} dY^g_{s+} + [M^X, M^Y]_t + \sum_{0\leq s<t} \Delta Y^g_s \Delta^+ X_s \\
	& = \int_{(0,t]}Y_{s-}dX^r_s + \int_{[0,t)}Y_{s} d X^g_{s+} + \int_{(0,t]} X_{s-}dY^r_s + \int_{[0,t)} X_{s+} dY^g_{s+}.
\end{align*}
To see that the finite variation part only moves on the set $\seq{XY = 0}$, it is sufficient to note that $\seq{XY \neq 0} = \seq{X\neq 0}\cap\seq{Y \neq 0}$.
\end{proof}

For simplicity, we present the following lemma for $C^1$-functions rather than bounded measurable, since unlike Nikeghbali \cite{N1} or Cheridito et al. \cite{CNP}, we do not attempt to solve the Skorokhod embedding problem for optional semimartingales of class-$(\Sigma)$ and we only include the following result to illustrate that class-$(\Sigma)$ is closed under this transform.

\vskip10pt
\bl\label{l3.2} Let $B= B^c + B^g$ be a left continuous increasing process such that $dB_+$ is carried on the set $\{X = 0\}$, then for any given $C^1$-function $f$
\begin{align*}
f(B_{t})X_t   & = f(0)X_0 + \int_{(0,t]} f(B_s)d(X^r_s - B^c_s) + \int_{[0,t)} f(B_{s+})d(X^g_s - B^g_s) + \int_{[0,t)} f(B_{s+})dB_{s+}.
\end{align*}
and $f(B_{+})dB_{+}$ is carried on the set $\{f(B)X = 0\}$.

\noindent In particular, if $X = M + A$ is of an optional semimartingale of class-$(\Sigma)$ and $B = A$ then
\begin{gather*}
f(A_{t})X_t = f(0)X_0 + \int_{(0,t]} f(A_{s})dM_s +  \int_{[0,t)} f(A_{s+})dA_{s+}
\end{gather*}
is an optional semimartingale of class-$(\Sigma)$.
\el 

\begin{proof}
We first note that the process $B = B^c  + B^g$ is left continuous and by an application of the It\^o formula in \cref{ito} we obtain
\begin{align*}
f(B_{t})X_t & = f(0)X_0 + \int_{(0,t]} f(B_{s})dX^r_s + \int_{[0,t)} f(B_{s})dX^g_s \\
		  & \quad + \int_{(0,t]} X_{s-}df(B)^r_s + \int_{[0,t)} X_{s}df(B)^g_{s+} + \sum_{s<t} (f(B_{s+})- f(B_s))\Delta A^+_s.
\end{align*}
By applying the It\^o formula to $f(B)$, we obtain
\begin{align*}
df(B_t) 
		& = \int_{(0,t]}f'(B_s)dB^c_s + \int_{[0,t)}f'(B_s)dB^g_{s+}  + \sum_{s<t} f(B_{s+}) - f(B_s) - f'(B_s)\Delta^+ B_s.
\end{align*}
The above implies that $df(B)_+$ is carried on the set $\seq{X=0}$ since $dB_+$ is carried on $\seq{X= 0}$. Therefore we have
\begin{align*}
 f(B_t)X_t & = f(0)X_0 + \int_{(0,t]} f(B_s)dX^r_s + \int_{[0,t)} f(B_s)dX^g_s + \sum_{s<t} (f(B_{s+})- f(B_s))\Delta^+ A_s\\
		   & = f(0)X_0 + \int_{(0,t]} f(B_s)dX^r_s + \int_{[0,t)} f(B_{s+})dX^g_s\\
   		   & = f(0)X_0 + \int_{(0,t]} f(B_s)d(X^r_s - B^c_s) + \int_{[0,t)} f(B_{s+})d(X^g_s - B^g_s) + \int_{[0,t)} f(B_{s+})dB_{s+}.
\end{align*}
It is sufficient to note that since $dB_+$ is carried on $\seq{X= 0 }$
\begin{gather*}
\int_{[0,t)} \I_\seq{X_sf(B_{s}) \neq 0}f(B_{s+})dB_{s+} = \int_{[0,t)} \I_\seq{X_s\neq 0}\I_\seq{f(B_{s})\neq 0}f(B_{s+})dB_{s+} = 0,
\end{gather*}
which concludes the proof.
\end{proof}

\subsection{The Madan-Royette-Yor formula}
From the works of Madan et al. \cite{MRY}, Profeta et al, \cite{PRY} and generalizations of Cheridito et al. \cite{CNP} we know that there is a deep connection between semimartingales of class-$(\Sigma)$ and their last passage time of zero. From this point onwards, given a optional semimartingale $X$ of class-$(\Sigma)$, the honest time $\tau := \sup\seq{s:X_s = 0}$ is assumed to be finite. Let us first record below Corollary 3.5 from \cite{CNP}.

\bp[Corollary 3.5 in \cite{CNP} or see Madan, Roynette and Yor \cite{MRY}]\label{PMRYF}
Let $K$ be a constant and M is a local martingale with no positive jumps such that $M^-$ is of class-$(D)$. Denote $g_K = \sup\{t \geq 0 : M_t= K\}$ then for every stopping time $T$,
\begin{gather}\label{MRYF}
(K-M_T)^+ = \mathbb{E}((K-M_\infty)^+ \I_{\{g_K\leq T\}} | \F_T).
\end{gather}
\ep
In the \cref{PMRYF}, the process $M$ is assumed to have no positive jumps, hence by \cref{l3.1} (ii) the process $X := (K-M)^+$ is a positive submartingale of class-($\Sigma$) and $g_K = \tau = \sup\seq{s:X_s = 0}$. Formula \eqref{MRYF} can then be rewritten into the form $\mathbb{E}(X_\infty \I_{\{\tau\leq T\}} | \F_T) = X_T$ which, with a slight abuse of terminology, we shall refer to as the {\it Madan-Royette-Yor type formula}. The goal now is to recover formulae of this type for optional semimartingale of class-$(\Sigma)$ and their last passage times of zero.

In the following, deviating from \cite{CNP, EOM, EOMT, N1}, we will work with both processes of class-$(DL)$ and class-$(D)$. To proceed, we introduce the time change processes
\begin{align*}
\tau_t   = \sup\seq{s <  t: X_s = 0}  \quad \textrm{and} \quad  \gamma_t = \inf\seq{s \geq  t: X_s = 0}.
\end{align*}
where by convention $\gamma_t$ takes the value infinity if the set is empty. We observe that $\tau_\infty:= \lim_{t\rightarrow \infty} \tau_t = \tau$ and $\seq{\tau_t < u} \subset \seq{t\leq \gamma_u} \subset \{\tau_t \leq u\}$ for $0\leq u <t \leq \infty$. Here the set inclusions are strict, since there could be trajectories such that $\gamma_u \geq t$ but $\tau_t = u$ and $X_{\tau_t} \neq 0$. However we have
\begin{align}\label{Xzero}
\seq{t\leq \gamma_u} \cap \{X_{\tau_t}  = 0\} =  \{\tau_t < u\} \cap \{X_{\tau_t}  = 0\} 
\end{align}
which follows from the fact that $\seq{t\leq \gamma_u} \cap \{X_{\tau_t}  = 0\} \cap \{\tau_t = u\}  = \emptyset$, since here $\gamma_u = u < t$. 
Similarly, we also consider the time change processes
\begin{align*}
g_t = \sup\seq{s \leq  t: X_s = 0} \quad \textrm{and} \quad k_t = \inf\seq{s >  t: X_s = 0}
\end{align*}
where we have $\seq{g_t \leq  u} = \seq{t< k_u}$ for $0\leq u < t$. It is important to point out that for $t= u$, the set $\seq{g_u \leq  u}$ is of probability one, but the set $\Omega \setminus \seq{u< k_u}$ could be of positive probability.

We derive below a balayage type formula for $X$ and $X_+$ which provide us with the martingale that will underpin later computations. Note that the right jumps $\Delta^+ X = \Delta A^g$ are only non-zero on the set $\seq{X= 0}$, this implies $X_{\gamma_t} = X_{k_t} = 0$ but the quantities $X_{\tau_t}$ and $X_{g_t}$ may or may not be zero. Also by definition $X_0 = 0$ and by convention $X_{0-} = 0$, hence $g_0 = 0$ and we can set $\tau_0 = 0$.

\bl\label{l3.4} Let $X$ be an optional semimartingale of class-$(\Sigma)$ then for $0\leq u \leq t < \infty$
\begin{align*}
X_{t+}\I_\seq{t <  k_u} & = M_{t\wedge k_{u}} +  A^c_{u} + A^g_{u+} - \Delta^+ A^g_u \I_{\{k_u = u\}}\\
X_t\I_\seq{t \leq \gamma_u} & = M_{t\wedge \gamma_{u}} +  A^c_{u} + A^g_{u}.
\end{align*}

\el 
\begin{proof}
For $0\leq u \leq t < \infty$, we have for $X_+$,
\begin{align*}
X_{t+}\I_\seq{t < 	k_u} 
	& = M_{t\wedge k_u} +  A^c_{t\wedge u} + (A^g_{+})_{t\wedge k_u}- X_{k_u+}\I_\seq{t\geq k_u}\\
	& = M_{t\wedge k_u} +  A^c_{t\wedge u} + (A^g_{+})_{t\wedge k_u}- \Delta A^g_{k_u+}\I_\seq{t\geq k_u}
\end{align*}
where in the last equality follows from the fact that $X_{k_u} = 0$ and $\Delta^+ X = \Delta A^g_+ = \Delta^+ A^g$. To further simply the above expression, we see that on the set $\seq{t\geq k_u}$ 
\begin{align*}
(A^g_{+})_{t\wedge k_u} - \Delta A^g_{k_u}\I_\seq{t\geq k_u} & = A^g_{k_u} \I_\seq{t\geq k_u}\\
							 & =  A^g_{u} \I_\seq{t\geq k_u} \I_\seq{u = k_u} + A^g_{k_u} \I_\seq{t\geq k_u} \I_\seq{u < k_u}\\
							 & =  A^g_{u+}\I_\seq{u = k_u}- \Delta^+A^g_{u} \I_\seq{u = k_u}+ A^g_{k_u} \I_\seq{t\geq k_u} \I_\seq{u < k_u}
\end{align*}
and, by using the fact that $A^g$ is left continuous, the third term above is given by
\begin{align*}
A^g_{k_u} \I_\seq{t\geq k_u} \I_\seq{u < k_u} & = A^g_{u+} \I_\seq{t\geq k_u} \I_\seq{u < k_u} \I_\seq{\Delta^+ A^g_u \neq 0 } + A^g_{u} \I_\seq{t\geq k_u} \I_\seq{u < k_u} \I_\seq{\Delta^+ A^g_u = 0 }\\
											  & = A^g_{u+}\I_\seq{t\geq k_u} \I_\seq{u < k_u}.
\end{align*}
On the complement $\seq{t < k_u}$, we have $(A^g_{+})_{t\wedge k_u} = A^g_{t+} = A^g_{u+}$, where the last equality comes from the fact that $u\leq t <k_u$ and $A^g_+$ does not increase on $\llb u, k_u\llb$. By combining the above computations, we obtain
\begin{align*}
& X_{t+}\I_\seq{g_t \leq 	u}  \\
& = M_{t\wedge k_u} +  A^c_{u} + A^g_{u+}\I_\seq{t<k_u} + A^g_{u+} \I_\seq{t\geq k_u} \I_\seq{u < k_u} +  A^g_{u+} \I_\seq{t\geq k_u} \I_\seq{u = k_u}  - \Delta^+ A^g_{u} \I_\seq{t\geq k_u} \I_\seq{u = k_u} \\
   						      						   & = M_{t\wedge k_u} +  A^c_{u} + A^g_{u+} - \Delta^+ A^g_{u}\I_\seq{u = k_u} 
\end{align*}

Similarly, by using the fact that $X_{\gamma_u} = 0$ on the set $\{\gamma_u <\infty\}$ we have
\begin{align*}
X_{t\wedge \gamma_u } = X_{t}\I_{\{t \leq \gamma_u\}} + X_{\gamma_u }\I_{\{t > \gamma_u\}} = X_{t}\I_{\{t \leq \gamma_u\}}.
\end{align*}
On the other hand, we deduce from the fact that $A^g$ is left continuous and $dA^g_+$ is carried on $\{X = 0\}$ 
\begin{align*}
X_{t\wedge \gamma_u } = M_{t\wedge \gamma_u }  + A^c_{t\wedge \gamma_u } + A^g_{t\wedge \gamma_u } = M_{t\wedge \gamma_u }  + A^c_{u} + A^g_{u} 
\end{align*}
which is a local martingale for on $t \in [u, \infty)$.
\end{proof}

The first observation we make is that $\I_\seq{u \leq \gamma_u}  = 1$ and $X_{u+}\I_{\{u=k_u\}} = 0$.  The equality $X_{u+}\I_{\{u=k_u\}} = 0$ follows from the fact that, on the set $\{u = k_u\}$ there exists a sequence of random times $(k_u^n)_{n\in \mathbb{N}}$ strictly greater than $u$ such that $X_{k_u^n} = 0$ and $\lim_{n\rightarrow \infty} k_u^n = u$, which implies that $\lim_{n\rightarrow \infty} X_{k_u^n} = X_{u+} =0$.
Therefore from \cref{l3.4}, we have for fixed $u\geq 0$, 
\begin{align}
X_{t+}\I_\seq{t < k_u}  - X_{u+}& = M_{t}^{k_{u}}   - M_{u} 	\qquad u\leq t,\label{bala0}\\
X_{t}\I_\seq{t\leq \gamma_u}  - X_{u}& = M_{t}^{\gamma_{u}}   - M_{u}	\qquad u \leq t.\label{bala}
\end{align}
By examining \eqref{bala0}, we note that if the local martingale $M^{k_u}-M^u$ is a true martingale on $[u,\infty)$ then one can take the conditional expectation with respect to $\F_u$ and eliminate the right hand side using optional sampling theorem (see for example Theorem 2.58 of \cite{HWY}). The second observation we make is that the integrability properties of $M^{k_{u}}_s - M_{s}^{u}$ for $s \in [u,\infty)$ can be derived from the integrability properties of $X$ or $X_+$. In view of this, one can, in the definition of class-$(\Sigma)$, restrict ourselves to optional semimartingales for which $M$ is a martingale, however the goal is to look for some sufficient conditions on the process $X$ or $X_+$. The assumption that $X$ and $X_+$ are of class-$(D)$ is likely too strong for problems on a finite horizon. For example, the Brownian motion is of class-$(\Sigma)$ but does not belong to class-$(D)$.

\bt\label{t3.1}
Let $X$ be a optional semimartingale of class-$(\Sigma)$ and $0\leq u \leq t < \infty$,\\
(i) if $X_+$ is of class-$(DL)$ then 
\begin{align*}
\mathbb{E}(X_{t+}\I_\seq{t < k_u}|\F_u) & = X_{u+}.
\end{align*}
(ii) If $X$ is of class-$(DL)$ then
\begin{align*}
\mathbb{E}(X_t\I_\seq{t \leq  \gamma_u}|\F_u) & = X_u, \quad  0\leq u\leq t < \infty,
\end{align*}

\begin{proof}
We prove only (ii) since the proof for (i) is similar. Given a fixed $u\geq 0$, we have from \eqref{bala} that for $u\leq t$
\begin{align*}
X_{t}\I_\seq{t \leq \gamma_u} - X_{u}\I_\seq{u \leq \gamma_u} = \int_{(0,t]} \I_{\{u < s\leq \gamma_u\}} dM_s =: M_t(u).
\end{align*}
The process $M(u)$ is a local martingale and $M_t(u) = (M^{\gamma_u}_t - M_u)\I_{\{t\geq u\}}$. To obtain the claim, it is enough to show that $M(u)$ is a martingale and then take the conditional expectation with respect to $\F_u$. To do that we make use the fact that a local martingale is a martingale if and only if it is of class-$(DL)$. To see that the local martingale $M(u)$ is of class-$(DL)$, we observe that for any $t \geq 0$ and any stopping time $T$
\begin{align*}
|M_{t\wedge T}(u)| & \leq |X_{t\wedge T}\I_\seq{t\wedge T \leq  \gamma_u} - X_{u}|\I_{\{t\wedge T \geq u\}}\\
					   & \leq |X_{t\wedge T}| + |X_{t\wedge u\wedge T}|.
\end{align*}
Since $X$ is of class-$(DL)$, we deduce that $M(u)$ of class-$(DL)$ and hence a martingale. 
\end{proof}
\et

\brem 
For positive optional submartingales, the Madan-Royette-Yor type formulae established above in \cref{t3.1} can be viewed as a special case of a general result on the {\it multiplicative system} associated with a positive optional submartingale. For interested readers, we refer to Lemma 3.10 in the recent work of Jeanblanc and Li \cite{JL} and the references within.
\erem

To establish the analogues formula at $t =\infty$, we point out that, since $\tau = \sup\seq{s:X_s = 0}$ is assumed to be finite and the measure $dA_+$ is carried on the set $\seq{X=0}$, the process $A$ is flat on $\rrb \tau, \infty \rrb$. Hence if $X_+$ converges almost surely to an integrable random variable $X_\infty$, then $X$ must also converge almost surely to $X_\infty$. In view of this, we simplify the problem and suppose below that $X_+$ is of class-$(D)$ and make use of limit results for c\`adl\`ag submartingales.

\bt\label{c3.1}
Let $X$ be a optional semimartingale of class-$(\Sigma)$ such that $X_+$ is of class-$(D)$ then\\ 
(i) for any finite stopping time $\sigma$, we have
\begin{align*}
\mathbb{E}(X_\infty\I_\seq{\tau \leq  \sigma}|\F_\sigma) &  = X_{\sigma+}
\end{align*}
(ii) and if $X_\tau =0$ then we have 
\begin{align*}
\mathbb{E}(X_\infty\I_\seq{\tau <  \sigma}|\F_\sigma) & = X_{\sigma}.
\end{align*}
\et
\begin{proof}
From \cref{l3.1}, we know that $X$ and hence $X_+$ can be written as the difference of two submartingales, that is $X_+ = (X^+)_+ - (X^-)_+$ and $|X_+| =  (X^+)_+ + (X^-)_+$. This implies that both $(X^+)_+$ and $(X^-)_+$ are right continuous positive submartingales of class-$(D)$ and there exists
$X_\infty := X^+_\infty-X^-_\infty \in L^1.$ Since $X^+$ is an optional submartingale of class-$(D)$, we have all stopping times $T$, the inequality $\mathbb{E}(X_\infty^+ \,|\,\F_T) \geq (X^+)_{T+} \geq X^+_{T}$. 

From this we deduce that $X^+$ is also of class-$(D)$ and $\lim_{t\rightarrow \infty} X^+_t =  X_\infty$. Similar arguments shows that $X^-$ is an optional submartingale of class-$(D)$.  From the Doob-Meyer-Mertens-Gal'\v cuk decomposition we can write $X^+ = m + a$ and $X^- = u+v$, where $m$ and $u$ are optional martingales and $a$ and $v$ are strongly predictable increasing process of integrable variation. From the decomposition $X = M + A$, we deduce that $M + A = m-u + a-v$. Then by using the fact that $X$ is of class-$(\Sigma)$, the process $A$ is left continuous and therefore strongly predictable, we see that $M-(m-u)= (a-v)-A$ is a c\`adl\`ag predictable local martingale of finite variation. This implies $M = m-u$ is an optional martingale and is therefore uniformly integrable. 

\noindent (i) From \cref{l3.4} or \eqref{bala0} and the fact that $\seq{g_t \leq  u} = \seq{t< k_u}$ for $0\leq u < t$, we deduce that for any finite stopping time $\sigma$,
\begin{gather*}
X_{\infty}\I_\seq{\tau\leq  \sigma}  - X_{\sigma+}= M_{k_{\sigma}}   - M_{\sigma}.
\end{gather*}
The result then follows from optional sampling theorem and the uniformly integrability of $M$.

\noindent (ii) We recall $\gamma_s = \inf\seq{u\geq s: X_u= 0}$ and observe that for any finite stopping time $\sigma$
\begin{align*}
X_{\gamma_{\sigma}} & = X_\infty\I_\seq{\tau< \sigma} + X_{\gamma_\sigma}\I_\seq{\tau\geq \sigma} = X_\infty\I_\seq{\tau< \sigma}
\end{align*}
where the second equality follows from the fact that $X_\tau = 0$ and
\bde
X_{\gamma_{\sigma}} \I_{\{\tau = \sigma\}} = X_{\infty} \I_{\{\tau = \sigma\}\cap\{X_\tau \neq 0\}} + X_{\tau} \I_{\{\tau = \sigma\}\cap\{X_\tau = 0\}} = X_{\infty} \I_{\{\tau = \sigma\}\cap\{X_\tau \neq 0\}} = 0.
\ede 
On the other hand, $X_{\gamma_\sigma} = M_{\gamma_\sigma} - A_{\sigma}$, since $A$ is left continuous and $dA_+$ is carried on $\seq{X= 0}$. The result again follows from applying optional sampling theorem to $M$.
\end{proof}

\brem 
As a check, one can apply \cref{c3.1} (i) and (ii) to $X := 1-\wt Z$, where $\wt Z$ is the Az\'ema supermartingale associated with a finite honest time $\tau$. From the fact that $\wt Z_\infty = Z_\infty = 0$ and $\wt Z_\tau = 1$, we recover $\P(\tau < \sigma \,|\, \F_\sigma) = 1-\wt Z_\sigma$ and $\mathbb{P}(\tau \leq  \sigma|\F_\sigma) = 1-Z_{\sigma}$.
Finally, one can also recover \cref{PMRYF} from \cref{c3.1} (i), by observing that for $X := (K-M)^+$ we have $X = X_+$ and $X_{k_\sigma} = 0$.
\erem

\subsection{Construction of finite honest times}\label{s3.2}
As an application of the results we have obtained on optional semimartingale of class-$(\Sigma)$, we show in \cref{max} a method to construct examples of optional submartingales of class-$(\Sigma)$, where both $A^c$ and the left continuous pure jump part $A^g$ are non-zero. This gives examples of finite honest times for which the representations obtain in \cref{t1} are non-trivial.

We recall that in continuous filtrations, all martingales are continuous, while for jumping filtration it was shown in Theorem 1 of Jacod and Skorokhod \cite{JS} that all martingales are almost surely of locally finite variation. Therefore by combining Theorem 1 of Jacod and Skorokhod \cite{JS} and \cref{l3.3}, one can produce non-trivial examples of positive optional submartingales of class-$(\Sigma)$ as defined in \cref{d3.2} by taking the products of known examples in the Brownian filtration (see Appendix A. in Mansuy and Yor \cite{MY}) and the Poisson filtrations (see Aksamit et al. \cite{ACJ}). 

\bd
A honest time is said to be of {\it type-c} (resp. {\it type-d}) if the martingale part of the Doob-Meyer-Mertens-Gal'\v cuk decomposition of $1-\wt Z$ is continuous (resp. locally of finite variation).
\ed 

\bp\label{max}
Let $\tau^c$ be a finite honest time of type-c and $\tau^d$ a finite honest time of type-d, then $\P(\tau^c\vee \tau^d <  t\,|\,\F_t) = \P(\tau^c <  t\,|\,\F_t) \P(\tau^d <  t\,|\,\F_t)$.
\ep

\begin{proof}
We see that both $X_t : = \P(\tau^c <  t\,|\,\F_t)$ and $Y_t:=\P(\tau^d <  t\,|\,\F_t)$ are positive optional submartingales of class-$(\Sigma)$. From \cref{honest} the random time $\tau^c\vee \tau^d$ is an finite honest time and from \cref{l3.3} the process $XY$ is a positive optional submartingale of class-$(\Sigma)$. We observe that
\begin{align*}
 \tau^c\vee \tau^d = \sup\seq{s: X_s = 0} \vee \sup\seq{s: Y_s = 0} = \sup\seq{s: X_sY_s = 0}
\end{align*}
and $(XY)_{\tau^c\vee \tau^d} = 0$ since $X_{\tau^c} = Y_{\tau^d} = 0$. The result then follows by an application of \cref{c3.1} (ii) to the process $XY$.
\end{proof}

The above result says that the Az\'ema optional submartingale associated with the maximum of two finite honest times can expressed as the product of the Az\'ema optional submartingale associated with each individual honest time. To the best of our knowledge, this type of representation has not previously appeared in the literature.

\bex\label{e3.1}
For examples of a honest time of type-c, we consider the following taken from Mansuy and Yor \cite{MY}. Let $B$ be a Brownian motion and
\begin{align*}
	 T_1 := \inf\seq{t: B_t = 1} \quad \mathrm{and}\quad \tau^c := \sup\seq{t\leq T :  B_t = 0}.
\end{align*}
The random time $\tau^c$ is a finite honest time and the Az\'ema's submartingale associated with $\tau^c$ is given by $X_t := \P(\tau^c <  t\,|\, \F^B_t)  = \P(\tau^c \leq t\,|\, \F^B_t) = M^c_t + A^c_t$ where
\begin{gather*}
M^c_t =  B^+_{t\wedge T_1}-\frac{1}{2}L^0_{t\wedge T_1}(B)  \quad \mathrm{and} \quad A^c_t = \frac{1}{2}L^0_{t\wedge T_1}(B).
\end{gather*}
Here the process $L^0(B)$ is the local time of the Brownian motion $B$ at zero. 

For a honest time of type-d, we consider the example in Proposition 4.12 of \cite{ACJ}. Let $J$ be a compound Poisson process with intensity $\mu$. Given $a \geq 0$, we set 
$$\tau^d := \sup\seq{t : \mu t - J_t \leq a}.$$
Then under certain conditions on the intensity and the distribution of the jump size, it is known that $\tau^d$ is a finite honest time. 

We denote by $\Psi(x)$ the ruin probability associated with the process $\mu t - J_t$, i.e., for every $x \geq 0$, $\Psi(x) := \P(t^x < \infty)$ with $t^x := \inf\seq{t : x + \mu t - J_t < 0}$.
  Then the Az\'ema submartingale and optional submartingale of $\tau^d$ admits the decomposition $Y_{t+} := \P(\tau^d \leq t\,|\,\F^J_t) = M^d_t + A^d_t$ and $Y_t := \P(\tau^d < t\,|\,\F^J_t) = M^d_t + A^d_{t-}$ where
\begin{align*}
M^d_t & = 1- (1 - \Psi(0))\sum_n \I_\seq{t\geq T_n} -	\Psi(\mu t - X_t - a)\I_\seq{\mu t - J_t\geq a} - \I_\seq{\mu t-J_t< a}\\
A^d_t & = (1 - \Psi(0))
\sum_n \I_\seq{t\geq T_n}.
\end{align*}
Here the martingale $M^d$ is of finite variation and $A^d$ is predictable pure jump process with jump times given by $(T_n)_{n\in \nn}$, where for $n> 1$,
\begin{align*}
T_1 := \inf\{t>0: \mu t - J_t = a\} \quad \mathrm{and} \quad T_n := \inf\{t > T_{n-1} : \mu t- J_t = a\}.
\end{align*}

Finally, we suppose that the Brownian motion $B$ and the compound Poisson process $J$ given above are independent of each other and we consider the joint filtration $\ff = (\F_t)_{t\geq 0}$ where $\F_t = \F^B_t\vee \F^J_t$. From \cref{max} we have $\wt Z_t = \P(\tau^c \vee \tau^d \geq t\, |\, \F_t) = (1-X_tY_t)$. To compute the multiplicative representation of $\wt Z$ obtained in \cref{t1}, it is sufficient to apply the It\^o formula to $1-XY$ to obtain $N$ and $\overline{N}$, which are given by 
$$N_t  = (1- X_tY_t) e^{\int^t_0 \frac{1}{2} Y_s dL^0_{s\wedge T_1}(B)} \quad \textrm{and} \quad \overline{N}_t = e^{\int^t_0 \frac{1}{2}Y_s dL^0_{s\wedge T_1}(B)}.$$
\eex

\noindent {\bf Acknowledgement:} The author wish to thank the anonymous referee for his/her readings and valuable advices on the writing of this paper.

\appendix
\setcounter{section}{0}
\renewcommand{\thesection}{\Alph{subsection}}
\section*{Appendix}
\addcontentsline{toc}{section}{Appendix}

\subsection{Theory of enlargement of filtration}

\bl[Lemma 4.3, Chapter IV of \cite{J2} or Lemma 1.53 and Proposition 1.54 in \cite{AJ}] \label{l1.2}\hfill\break
Given a random time $\tau$, we have that\hfill\break 
(i) the sets $\seq{\wt Z_- = 1}$ and $\seq{\wt Z = 1}$ are the largest predictable and optional set contained in the stochastic interval $\llb 0,\tau\rrb$,\hfill\break 
(ii) the stochastic interval $\llb 0,\tau\rrb$ is contained in the sets $\seq{\wt Z > 0}$ and $\seq{\wt Z_- > 0}$,\\
(iii) the stochastic interval $\llb 0,\tau\,\llb$ is contained in the sets $\seq{Z > 0}$ and $\seq{\,^pZ > 0}$.
\el 

\bp[Proposition 5.1, Chapter V of \cite{J2} or Theorem 5.8 in \cite{AJ}] \label{p1.1}\hfill\break
The following are equivalent;\hfill\break 
(i) a random time $\tau$ is a finite honest time.\\
(ii) $\tau = \sup \seq{s: \wt Z_s = 1}$, i.e. it is the end of an optional set.\\
(iii) $H^o_t = H^o_{t\wedge \tau}$ for $t\geq 0$.\\
(iv) $\wt Z_\tau = 1$.
\ep 
\bl\label{l1.3}
Let $\tau$ be a finite honest time then $dH^o$ is carried on $\seq{\wt Z=1}$.
\begin{proof}
It enough to note that $\mathbb{E}(\int_{[0,\infty)} \I_\seq{\wt Z_s< 1} dH^o_s) = \P(\wt Z_\tau< 1) = 0$
\end{proof}
\el

\subsection{Stochastic calculus for optional semimartingales}

\bd\label{d1.1}
A stochastic process $X$ is said to be an optional (super or sub)martingale if it is a (super or sub)martingale and (i) $X$ is an optional process, (ii) for any stopping time $T$, $X_T\I_\seq{T<\infty}$ is integrable, (iii) there exists an integrable random variable $\zeta$ such that for any stopping time $T$, $X_T = \mathbb{E}(\zeta|\F_T)$ ($X_T \geq \mathbb{E}(\zeta|\F_T)$ or $X_T \leq \mathbb{E}(\zeta|\F_T)$) a.s. on the set $\seq{T<\infty}$. 
\ed 

\bd\label{d1.3}
A l\`agl\`ad stochastic process $X$ is said to be strongly predictable if $X$ is predictable and the right limit $X_+$ is optional.
\ed 

\bd
A stochastic process $X$ is said to be an optional semimartingale if $X$ can be written as 
\bde
X = X_0 + M + A, \quad M_0 =0,\quad  A_0 = 0,
\ede 
where $M$ is a local martingale and $A$ is an l\`agl\`ad adapted process of finite variation.
\ed

\bt[It\^o formula. Theorem 8.2 \cite{G1}]\label{ito} Let $X = (X^1,\dots, X^k)$ be an optional semimartingale and $X^k = X^k_0 + M^k + A^k$ for $k= 1,\dots,n$. Let $F(x) = F(x_1,\dots, x_n)$ be a continuously twice differentiable function on $\rr^n$ then for $t\in \rr_+$,
\begin{align*}
F(X_t) &= F(X_0) + \sum_{k=1}^n \int_{(0,t]} D^kF(X_{s-})d(A^{k,r} + M^{k,r})_s + \frac{1}{2} \sum_{k,l=1}^n\int_{(0,t]} D^kD^lF(X_{s-}) d\left<M^{k,c},M^{l,c}\right>_s\\
	   &\quad + \sum_{0<s\leq t} \left[F(X_s) - F(X_{s-}) - \sum^n_{k=1} D^kF(X_{s-})\Delta X^k_s \right] +\sum_{k=1}^n \int_{[0,t)} D^kF(X_s)dA^{k,g}_{s+} \\
	   & \quad + \sum_{0\leq s< t} \left[F(X_{s+}) - F(X_{s}) - \sum^n_{k=1} D^kF(X_{s})\Delta^+ X^k_s \right]
\end{align*}
where $D^k$ is the partial derivative with respect to the $k$-th coordinate and $M^r = M^c + M^d$.
\et

\bl[Tanaka formula. Lemma 5.7 \cite{GIOQ}] \label{tanaka} Let $X$ be a (real-valued) optional semimartingale
with decomposition $X = X_0 + M + A^r + A^g$. Let $f : \rr \rightarrow \rr$ be a convex
function. Then $f(X)$ is an optional semimartingale. Moreover, denoting by $f'$ the left-hand
derivative of $f$, then we have
\begin{align*}
f(X_t) &= f(X_0) + 
\int_{(0,t]} f'(X_{s-})d(A^r_s + M_s) + \int_{[0,t)}f'(X_s)dA^g_{s+}\\
& \quad  + \sum_{0<s\leq t} f(X_s) - f(X_{s-}) - f'(X_{s-})\Delta X_s  + \sum_{0\leq s < t}  f(X_{s+}) - f(X_s) - f'(X_s)\Delta^+X_s + C^f_t
\end{align*}
where $C^f$ is a continuous increasing process.
\el

\bt[Dol\'eans-Dade-exponential. Theorem 5.1 \cite{G3}] \label{stochexp}
Let $X$ be an optional semimartingale. There exists a unique (to within indistinguishably) optional semimartingale $S$ such that 
\bde
S_t = S_0 + \int_{(0,t]} {S_{s-}}dX^r_s + \int_{[0,t)} S_s dX^g_{s+}
\ede 
The process $S$ is given by the formula
\bde
S = S_0\exp\left\{ X - \frac{1}{2}\left< X^c,X^c\right> \right\}\prod_{0<s\leq \cdot} (1+\Delta X_s)e^{-\Delta X_s}\prod_{0<s< \cdot} (1+\Delta^+ X_s)e^{-\Delta^+ X_s}
\ede 
and is termed the optional stochastic exponential $X$ which we shall denote by $\mathcal{E}(X)$.
\et 

\bt[Doob-Meyer-Mertens-Gal'\v cuk decomposition, \cite{G2} \cite{MJ}]\label{dm1}
An optional supermartingale $X$ admits a decomposition $X = M-A$, where $M$ is a (local) optional martingale and $A$ is an increasing strongly predictable (locally) integrable process with $A_0 = 0$ if and only if $X$ belongs to the class-(D) (class-(DL)). This decomposition is unique to within indistinguishably.
\et

To our knowledge, under the usual conditions, the {\it local time} of optional semimartingales has not been studied. 
By following similar arguments to section 6, Chapter IX of \cite{HWY}, we define the local time $L^a(X)$ of an optional semimartingale $X$ at $a \in \rr$ and show that $dL^a(X)$ is carried on $\seq{X=a}$.

\bl\label{tanaka2}
Let $X$ be an optional semimartingale and $a\in \rr$. Then 
\begin{align*}
(X_t-a)^+ &= (X_0-a)^+ + 
\int_{(0,t]} \I_\seq{X_{s-} > a}d(A^r_s + M_s) + \int_{[0,t)} \I_\seq{X_{s} > a}\, dA^g_{s+} \\
 & \quad + \sum_{0< s \leq  t}  \I_\seq{X_{s-} > a }(X_{s}-a)^- + \sum_{0< s \leq  t} \I_\seq{X_{s-} \leq a }(X_{s}-a)^+ \\
 & \quad + \sum_{0\leq s < t}  \I_\seq{X_s > a }(X_{s+}-a)^- + \sum_{0\leq s < t}   \I_\seq{X_s \leq a }(X_{s+}-a)^+  + \frac{1}{2}L^a_t(X),\\
(X_t-a)^- &= (X_0-a)^- + 
\int_{(0,t]} \I_\seq{X_{s-} \leq a}d(A^r_s + M_s) + \int_{[0,t)} \I_\seq{X_{s} \leq a}\, dA^g_{s+} \\
 & \quad + \sum_{0< s \leq  t}  \I_\seq{X_{s-} > a }(X_{s}-a)^- + \sum_{0< s \leq  t} \I_\seq{X_{s-} \leq a }(X_{s}-a)^+\\
 & \quad  + \sum_{0\leq s < t}  \I_\seq{X_s > a }(X_{s+}-a)^- + \sum_{0\leq s < t}   \I_\seq{X_s \leq a }(X_{s+}-a)^+  + \frac{1}{2}L^a_t(X)
 \end{align*}
where $L^a(X)$ is a continuous adapted increasing process with $L^a_0(X) = 0$. The process $L^a(X)$ is called the local time of $X$ at $a\in\rr$.
\el

\begin{proof}
By applying the Tanaka formula in \cref{tanaka} to $f(x) = (x-a)^+$ and $g(x) = (x-a)^-$ and taking the difference, we obtain
\begin{align*}
(X_t-a)^+ - (X_t-a)^- = (X_0-a)^+ - (X_0-a)^- + \int_{(0,t]} dX^r_s + \int_{[0,t)} dX^g_{s+} + C^f_t-C^g_t.
\end{align*}
This gives $C^f = C^g$ and we denote them by $\frac{1}{2}L^a(X)$.
\end{proof}
\bt\label{tanaka3}
Let $X$ be an optional semimartingale and $a \in \rr$ then the support of $L^a(X)$ is contained in $\seq{X = a}$.
\et 
\begin{proof}
The proof is similar to that of Theorem 9.44 of \cite{HWY} in the c\`adl\`ag case. We start by supposing $0<S\leq T$ and $\llb S,T\llb \subset \seq{X < 0}$. This implies that $\rrb S,T\rrb \subset \seq{X_- \leq 0}$ and $\llb S,T\llb \subset \seq{X_+ \leq 0}$. By applying the Tanaka formula in \cref{tanaka},
\begin{align*}
(X_T-a)^+ - (X_S-a)^+  &= \sum_{S< s \leq  T} \I_\seq{X_{s-} \leq a }(X_{s}-a)^+ + \frac{1}{2}L^a_T(X) - \frac{1}{2}L^a_S(X)\\
	  &= (X_{T}-a)^+ + \frac{1}{2}L^a_T(X) - \frac{1}{2}L^a_S(X)
	  \end{align*}
which implies that $L^a_T(X)  = L^a_S(X)$. Let $r$ be a rational number and set 
\begin{align*}
S(r) &=
\begin{cases}
r & X_{r} < a\\
\infty & X_{r} \geq a
\end{cases}\\
T(r) &= \inf \seq{t>S(r): X_{t} \geq a}\\
H 	 &= \bigcup_{r>0} \rrb S(r), T(r)\llb 
\end{align*}
For each $\omega \in \Omega$, the section $H(\omega)$ is the interior of the $\seq{t:X_t(\omega)<0}$. We see that the process $L^a(X)$ does not increase in the interior of $\seq{t:X_t(\omega)<0}$ and by similar arguments $L^a(X)$ does not increase on the interior of $\seq{t:X_t(\omega)>0}$. We conclude by noticing that the set $\seq{t:X_t(\omega)\neq 0}$ differs from its interior by a countable set, since the boundary set is contained in the set of jumps of the optional semimartingale $X$ which for each $\omega \in \Omega$ is countable (see Theorem 1.14 \cite{G1}). 
\end{proof}

\subsection{Proof of Lemma \ref{srl}} \label{4.3}
\begin{proof}
We first show that the pair $(z^*, a^*)$ where $z^* = y+ a^*$ and $a^* = (\overline{-y} \vee 0) =\overline{-y}= - \underline{y}$ is a solution to the Skorokhod reflection problem. To see that $a^*(t)=\sup_{s \le t}-y(s)$ is indeed a solution. Condition (i) is clearly satisfied since for all $t\geq 0$ 
$$z^*(t)=y(t) - \inf_{s \le t}y(s) \ge 0$$ 
and condition (ii) is satisfied since $y(0)=0$ and $y$ has continuous running infimum. 

To show that condition (iii) hold, we use the fact that $a^*$ is a continuous non-decreasing process and the number of left and right jumps of $y$ is at most countable (see Theorem 5.64 in \cite{TBB}) to obtain
\begin{align*}
\int_{[0,\infty)}  \I_\seq{-y(s) < a^*(s)}  da^*(s)  = \int_{[0,\infty)}  \I_\seq{\Delta^+ y(s)\Delta y(s) = 0 } \I_\seq{-y(s) < a^*(s)}  da^*(s).
\end{align*}
To see that the above integral is zero, it is sufficient to notice that if $s_0$ is a point of continuity of $-y$ and is such that $-y(s_0) < a^*(s_0)$ then there exists $\epsilon(s_0)>0$ such that $-y(s) < a^*(s)$ for all $s\in [s_0-\epsilon(s_0), s_0+\epsilon(s_0)]$. This implies that $a^*(s) = a^*(s_0)$ for all $s \in [s_0-\epsilon(s_0), s_0+\epsilon(s_0)]$ and $\seq{s:\Delta^+ y(s)\Delta y(s) = 0 } \cap \seq{s:-y(s) < a^*(s)} \subseteq S(a^*)^c$.

To prove uniqueness, we first show that $a^*$ is the smallest increasing function satisfying the condition that $z^*=y+a^* \ge 0$. To see this, for any $s\leq t$, we note that any solution $(z,a)$ of the Skorokhod reflection problem must satisfy $a(t) \ge 0$ and
\begin{align*}
a(t) &\ge a(s) = z(s) - y(s) \ge -y(s),
\end{align*}
from which we deduce the inequality 
$$a(t) \ge \sup_{s \le t}(-y(s) \vee 0)=\sup_{s \le t}-y(s)=a^*(t).$$ 

To this end, suppose one is given another solution $(a,z)$ distinct from $(a^*, z^*)$, then there exists $t> 0$ such that $a(t) > a^*(t)$. For such $t>0$, we consider the time change process $$g(t)=\sup\{s\leq t:a^*(s)=a(s)\},$$
and note that $g(t)$ is well defined since $a$ and $a^*$ are both continuity and starts at zero. The interval $(g(t),t]$ must be non-empty since $a(0) = 0 = a^*(0)$, $a(t) > a^*(t)$ and both $a$ and $a^*$ are continuous (use intermediate value theorem). Moreover, by continuity, we have $a\left(g(t)\right)=a^*\left(g(t)\right)$ and for $s\in \left(g(t),t\right]$ we must have $z(s)>z^*(s)\ge 0$ because $a(s)>a^*(s).$  However this is a contradiction as $a(s)$ cannot increase on $\left(g(t),t\right]$ since $\left(g(t),t\right] \subset \seq{s:z(s)>0}$. Hence $(a^*, z^*)$ is the unique solution.
\end{proof}

\end{document}